\documentclass[12pt]{article}
\usepackage{amsmath}
\usepackage{amssymb}
\usepackage{amsfonts}
\usepackage{dsfont}
\usepackage{bm}
\usepackage{color,soul}
\usepackage{subfig}
\usepackage{amsthm}
\usepackage[final]{pdfpages}
\usepackage{multirow}
\usepackage{algorithm}
\usepackage[noend]{algpseudocode}

\newcommand{\blind}{1}
\usepackage{natbib}
\newtheorem{theorem}{Theorem}
\newtheorem{corollary}{Corollary}

\newtheorem{lemma}{Lemma}
\newtheorem{remark}{\bf Remark}

\newtheorem{assumption}{Assumption}
\usepackage{enumerate}

\usepackage{array}
\usepackage{mathtools}

\DeclarePairedDelimiter\abs{\lvert}{\rvert}%
\DeclarePairedDelimiter\norm{\lVert}{\rVert}%

\makeatletter
\let\oldabs\abs
\def\abs{\@ifstar{\oldabs}{\oldabs*}}
\let\oldnorm\norm
\def\norm{\@ifstar{\oldnorm}{\oldnorm*}}
\makeatother

\newcommand{\mD}{\mathcal{D}}

\newcommand{\be}{\boldsymbol{e}}

\newcommand{\bI}{\boldsymbol{I}}
\newcommand{\bJ}{\boldsymbol{J}}

\newcommand{\mL}{\mathcal{L}}

\newcommand{\bS}{\boldsymbol{S}}

\newcommand{\bU}{\boldsymbol{U}}
\newcommand{\bu}{\boldsymbol{u}}

\newcommand{\bv}{\boldsymbol{v}}

\newcommand{\bw}{\boldsymbol{w}}
\newcommand{\bX}{\boldsymbol{X}}
\newcommand{\bx}{\boldsymbol{x}}
\newcommand{\bY}{\boldsymbol{Y}}
\newcommand{\by}{\boldsymbol{y}}
\newcommand{\bz}{\boldsymbol{z}}

\newcommand{\bbeta}{\boldsymbol{\beta}}

\newcommand{\bgamma}{\boldsymbol{\gamma}}

\newcommand{\bSigma}{\boldsymbol{\Sigma}}

\newcommand{\bTheta}{\boldsymbol{\Theta}}

\title{
Statistical Inference in High-dimensional Generalized Linear Models with Streaming Data}

\begin{document}

\def\spacingset#1{\renewcommand{\baselinestretch}%
	{#1}\small\normalsize} \spacingset{1}


\if1\blind
{
	\title{\bf  Statistical Inference in High-dimensional Generalized Linear Models with Streaming Data}
	\author{
		Lan Luo\thanks{Equal contribution. Department of Statistics and Actuarial Science, University of Iowa, USA},\hspace{0.2 cm}
		Ruijian Han\thanks{Equal contribution. Department of Statistics, The Chinese University of Hong Kong, Hong Kong SAR, China},\hspace{0.2cm}
		Yuanyuan Lin\thanks{Department of Statistics, The Chinese University of Hong Kong, Hong Kong SAR, China},\hspace{0,2 cm}
		and   
		Jian Huang\thanks{ Department of Statistics and Actuarial Science, University of Iowa, Iowa, USA}
	}
	\date{\today}	
	\maketitle
} \fi

\if0\blind
{
	\bigskip
	\bigskip
	\bigskip
	\begin{center}
		{\LARGE\bf Statistical Inference in High-dimensional Generalized Linear Models with Streaming Data}
	\end{center}
	\medskip
} \fi

\begin{abstract}
	
	In this paper we develop an online statistical inference approach for high-dimensional generalized linear models with streaming data for real-time estimation and inference.  We propose an online debiased lasso (ODL)  method to accommodate the special structure of streaming data. ODL differs from offline debiased lasso in two important aspects. First, in computing the estimate at the current stage, it only uses
	summary statistics of the historical data. Second, in addition to debiasing an online lasso estimator, ODL corrects an approximation error term arising from nonlinear online updating with streaming data.
	We show that  the proposed online debiased estimators for the GLMs are consistent and asymptotically normal. This result provides a theoretical basis for carrying out real-time interim statistical inference with streaming data.  Extensive numerical experiments are conducted to evaluate the performance of the proposed ODL method. These experiments demonstrate the effectiveness of our algorithm and support the theoretical results. A streaming dataset from the
	National Automotive Sampling System-Crashworthiness Data System is analyzed to illustrate the application of the proposed method.
\end{abstract}

\noindent%
{\it Keywords:} Adaptive tuning, confidence interval, high-dimensional data, online debiased lasso, online error correction.

\newpage
\spacingset{1.5} 

\section{Introduction}
Due to the explosive growth of data from non-traditional sources such as sensor networks, security logs and web applications, streaming data is becoming a core component in big data analyses. By streaming data, we refer to data that is generated continuously over time, typically in high volumes and at high velocity. It includes a wide variety of data types such as log files generated by mobile or web applications, ecommerce purchases, information from social networks, and financial trading floors. To reduce the demand on computing memory and achieve real-time processing, the nature of streaming data calls for the development of algorithms which require only one pass over the data. 
Furthermore, a lot of data streams are high-dimensional in nature such as genomic data for explaining the variability in the phenotypes of an organism and its genetic profile~\citep{Genetic2017,Bioinform2017}, and neuroimaging data for predicting neural activity patterns given various stimulus~\citep{encode2018,encoding_correlated_2020}. At an initial stage of data collection, it is routine to encounter a case where the number of predictors exceeds the number of observations. Moreover, even if sample size grows along the data streams, the candidate set of predictors may contain a large portion of sparse redundancy features. Therefore, to improve interpretability of the results, we need to differentiate predictors in terms of their signal strength and one of the classical approaches is the penalized regression method.
In this paper, we focus on a high-dimensional regression setting in generalized linear models (GLM) with non-Gaussian outcomes. Our goal is to develop a real-time estimation and inference procedure that is highly scalable with respect to fast growing data volumes, but with no loss of efficiency in statistical inference in existence of a large number of features.

Streaming data processing essentially falls into the field of online learning. This line of research may be dated back five decades or so when~\citet{Robbins1951} proposed a stochastic approximation algorithm that laid a foundation for the popular stochastic gradient descent (SGD) algorithm~\citep{Sakrison1965}. Later on, the SGD algorithm and its variants have been extensively studied for online estimation and prediction~\citep{Toulis2017}, but the work of developing online statistical inference remains unexplored. A recent paper by~\citet{Fang2019} proposed a perturbation-based resampling method to construct confidence intervals for SGD, but it does not achieve desirable statistical efficiency and may produce misleading inference in high-dimensional settings. In addition to the SGD types of recursive algorithms, several online updating methods have been proposed to specifically perform sequential updating of regression coefficient estimators, including the online least squares estimator (OLSE) for the linear model, the cumulative estimating equation (CEE) estimator, the cumulatively updated estimating equation (CUEE) estimator by~\citet{Schifano2016CUEE} and the renewable estimator by~\citet{Luo2020} for nonlinear models.

Most of the aforementioned online algorithms are developed under the low-dimensional settings where the number of features is far less than the total sample size. However, as pointed out above, a prominent concern in high-dimensional streaming data analysis is that only a subset of the variables have nonzero coefficients. Besides small sample size issue at the early stage of data collection, processing such data stream without properly accounting for the sparsity in feature set may introduce significant bias and invalid statistical inference. It is worth noting that even if the cumulative sample size exceeds the number of features as time goes by, traditional estimation methods in low-dimensional settings such as maximum likelihood estimation (MLE) may still incur large bias especially in GLMs~\citep{2019-highdim-MLE}. Therefore, current state-of-art online learning algorithms in low-dimensional settings may be insufficient for processing high-dimensional data streams.

In traditional offline settings,
many methods have been developed for analyzing high-dimensional static data. Most of the work on variable selection in high-dimensional regression problems is along the line of lasso~\citep{tibshirani1996regression}, the Smoothly Clipped Absolute Deviation (SCAD) penalty~\citep{fan2001variable}, the adaptive lasso~\citep{zou2006adaptive}, and the minimax convex penalty (MCP)~\citep{zhang2010nearly}.
However, variable selection methods focus on point estimation without quantifying the uncertainty in estimates. Later on, statistical inference problems in high-dimensional settings, including interval estimation and hypothesis testing, have attracted much attention
since the pioneering works of~\cite{Zhang_delasso_2014},~\cite{vandegeer2014},~\cite{javanmard2014confidence},
~\cite{belloni2015uniform}, among others. Other important methods on inference for high-dimensional linear models include~\cite{cai2017confidence, tony2020semisupervised},~\cite{belloni2019valid}, etc.~\cite{vandegeer2014} extended the de-sparsified lasso to high-dimensional GLMs.~\cite{ning2017general} proposed to construct confidence intervals for high-dimensional M-estimators based on decorrelated score statistic. Recently, a novel splitting and smoothing inference approach for high-dimensional GLMs was proposed by~\cite{fei2021estimation}.

While significant progress has been made on statistical inference for high-dimensional regression problems under the traditional offline settings, variable selection and statistical inference for high-dimensional models with streaming data is still at its infancy stage. ~\citet{Sun2018ANF} introduced a systematic framework for online variable selection based on some popular offline methods such as MCP, elastic net~\citep{Zou2005Elastic}. But their focus is not on statistical inference. Different from this work, there are some existing methods considering the problem of inference. For example, ~\citet{Deshpande2019OnlineDF} proposed a class of online estimators for high-dimensional auto-regressive models. One of the most relevant works is a novel inference procedure in GLMs based on recursive online-score estimation~\citep{shi2020statistical}. However, in both works, the entire dataset is assumed to be available at an initial stage for computing an initial estimator, e.g. the lasso estimator; thereafter, a recursively forward bias correction procedure is conducted along sequentially arrived data points. However, the availability of the entire dataset at an initial stage is not a natural setup in online learning. To address this issue,~\citet{han2021online} proposed an online debiased lasso method for statistical inference in high-dimensional linear models with streaming data.

Unlike the case of high-dimensional linear models where the loss function depends on data only through sufficient statistics~\citep{han2021online}, parameters and data are not linearly separable in GLMs. Motivated by the renewable estimation method in low-dimensional GLMs~\citep{Luo2020}, we start off by taking a first-order Taylor expansion on the quadratic loss function to bypass the need of historical individual-level data. The key idea centers around using ``approximate summary statistics" resulting from Taylor expansions. However, this is not a trivial extension of the methods developed under low-dimensional settings. In high-dimensional settings where predictors are spuriously correlated, a data-splitting strategy is typically used for decorrelation where variable selection and estimation are conducted using two different sub-datasets~\citep{ning2017general,shi2020statistical,fei2021estimation}. A prominent concern of using such approximate summary statistics that involve previous estimates is that it may incur dependency in the corresponding estimating equation. Theoretically speaking, the dependency among recursively updated estimators poses extra technical challenge in establishing the non-asymptotic error bound. In our proposed online method for real-time confidence interval construction, we aim to address the following questions: (i) what types of approximate summary statistics to be stored to carry out an online debiasing procedure? (ii) will the error accumulate along the updating steps if we use the approximate summary statistics? (iii) will the online debiasing procedure maintain similar oracle properties to its offline counterpart? and (iv) how to choose the tuning parameter adaptively in an online setting where cross-validation that relies on splitting the entire dataset is not feasible.

The focus of this paper is to develop an online debiased lasso (ODL) estimator in high-dimensional generalized linear models with streaming datasets for real-time estimation and inference. Our new contributions include: (i) we propose a two-stage online estimation and debiasing framework that aligns with streaming data collection scheme; (ii)  ODL accounts for
sparsity feature in a candidate set of predictors and provides valid statistical inference results; and (iii)
ODLs for the GLMs are shown to be consistent and asymptotically normal. This result provides a theoretical basis for carrying out real-time interim statistical inference with streaming data.
ODL is inspired by the offline debiased lasso method \citep{Zhang_delasso_2014, vandegeer2014, javanmard2014confidence},
 however, it differs from the offline debiased lasso in two important aspects. First, in computing the estimate at the current stage, it only uses summary statistics of the historical data. Second, in addition to debiasing an online lasso estimator, ODL corrects an approximation error term arising from
online updating with streaming data. This correction is crucial to guarantee the asymptotic normality of the ODL estimator.


This paper is organized as follows. Section~\ref{sec:method} introduces the model formulation followed by our proposed online two-stage debiasing method to process high-dimensional streaming data. Section~\ref{sec: asymptotic} includes some large sample properties concerning the theoretical guarantees for our proposed method. Simulation experiments are given in Section~\ref{sec:sim} to evaluate the performance of our proposed method in comparison to both MLE and offline debiased estimator. We illustrate the proposed ODL method and apply it
to analyze a real data example in Section~\ref{sec:data}. Finally, we make some concluding remarks in Section~\ref{sec:discussion}. All technical proofs are deferred to the Appendix.

\subsection{Notation}
For a matrix $ \bX \in \mathbb{R}^{n\times p}$, we let $ \bX_{i\cdot}, \bX_{\cdot j} $ and $ \bX_{ij} $ denote the $ i$-th row, $ j$-th column and $ (i,j)$-element of matrix $ \bX. $
$ \bX_{i,-j} $ is a sub-vector of $ \bX_{i\cdot} $ with the $ j $-th element deleted and $ \bX_{-i,-j} $ is a sub-matrix of $ \bX $ with the $ i$-th rows and the $ j $-th column deleted while other elements remain unchanged. For a vector $ \bx \in \mathbb{R}^m $, we define its $ \ell_q $-norm as $ \lVert \bx \lVert_q = (\sum_{i=1}^p x_i^q)^{\frac{1}{q}}$. For a sequence of random variables $ \{\xi_n\}_{n\in\mathbb{N}} $ and a corresponding sequence of constants $ \{a_n\}_{n\in\mathbb{N}} $. We say that $ \xi_n = \mathcal{O}_p(a_n) $ if for any $ \epsilon > 0 $, there exist two finite numbers $ M,N > 0 $ such that $ P(|\xi_n/a_n|>M) < \epsilon $ for any $ n > N. $ Generally speaking, $ \xi_n = \mathcal{O}_p(a_n) $ denotes $ \xi_n/a_n $ is stochastically bounded. $ \xi_n = o_p(a_n) $ means that $ \xi_n/a_n $ converges to zero in probability. With the consideration of the streaming data, we use $ \bX^{(j)} $ and $ \bY^{(j)} $  to stand for $ \bX$ and $\by $, namely explanatory variable and explained variables, arriving in $ j$-th batch respectively. In addition, $ \bX^{(j)}_\star $ and $ \bY^{(j)}_\star $ (with star index) are the cumulative variables of $ \bX^{(j)}  $ and $ \bY^{(j)} $. For example, $ \bX^{(j)}_\star = ((\bX^{(1)})^\top, \ldots, (\bX^{(j)})^\top)^\top. $

\section{Methodology}\label{sec:method}
In this section,  we describe the proposed estimation method with streaming data, including the online lasso estimation and online debiased lasso estimation. We also provide an adaptive tuning method to select the regularization parameter. A rundown of our algorithm is summarized at the end of this section. 

Consider up to a time point $b\geq2$, there is a total of $N_b$ samples arriving in a sequence of $b$ data batches, denoted by $\mD^\star_{b} = \{\mD_{1},\dots,\mD_{b} \}$, and each contains $n_j=|\mD_{j}|$ samples, $j=1,\dots,b$.
Assume each observation $y_i$ is independently from an exponential dispersion model with density function
\begin{align*}
f(y_i\mid\bx_i,\bbeta^0,\phi) = a(y_i;\phi) \exp\left\{-\frac{1}{2\phi} d (y_i;\mu_i) \right\}, \ i \in\mD^\star_{b},
\end{align*}
where $d(y_i;\mu_i)$ is the unit deviance function involving the mean parameter $\mu_i=\mathbb{E}(y_i\mid\bx_i)$, and $a(\cdot)$ is a normalizing factor depending only on the dispersion parameter $\phi > 0$. The systematic component of a GLM takes the form $\mu_i = \mathbb{E}(y_i\mid\bx_i) = g(\bx_i^\top\bbeta^0)$ where $g(\cdot)$ is a known link function. The underlying regression coefficient $\bm{\beta}^0\in\mathbb{R}^p$ is of our interest,
which is assumed to be sparse with $s_0$ nonzero elements.
Specifically, we let $S_0:=\{r:\beta^0_{r}\neq 0\}$ be the active set of variables and its cardinality is $s_0$.  Our main goal is to conduct a point-wise statistical inference for the components of the parameter vector $\beta_{r}^0 \ (r=1,\dots,p)$ upon the arrival of every new data batch $\mD_b,\ b = 1,2,\dots$.
The log-likelihood function for the cumulative dataset $\mD^\star_{b}$ is
\[
\ell(\bbeta,\phi;\mD^\star_{b}) =\frac{1}{N_b}\sum_{i\in\mD^\star_{b}} \log f(y_i\mid\bx_i,\bbeta,\phi) = \frac{1}{N_b}\sum_{i\in\mD^\star_{b}}\log a(y_i;\phi) - \frac{1}{2N_b\phi} \sum_{i\in\mD^\star_{b}} d(y_i;g(\bx_i^\top\bbeta)).
\]
Based on $\mD_b^\star$, the standard offline lasso
 estimator  is defined as
\begin{equation}\label{lasso_offline}
\bar{\bbeta}^{(b)} = \underset{\bbeta\in\mathbb{R}^p}{\arg\min}
\left\{\frac{1}{2N_b}\sum_{i\in\mD^\star_{b}}d\left(y_i; g(\bx_i^\top\bbeta)\right)  + \lambda_b\|\bbeta \|_1\right\},
\end{equation}
where $ N_b = \sum_{j=1}^{b}n_j$ is the cumulative sample size and $\lambda_b$ is the regularization parameter 
However, as discussed in \cite{Luo2020} and \cite{han2021online}, the optimization of \eqref{lasso_offline} requires the historical data $ \mD^{\star}_{b-1} $ which is not observed in the streaming data.
So the standard lasso estimator cannot be computed.
Therefore,  new online method to find the lasso estimator and the debiased lasso tailored for streaming data is desired. The latter is for the purpose of providing a statistical inference. For the sake of clarity, we refer to the lasso estimator in \eqref{lasso_offline} as the {\it offline lasso estimator}.
Sections \ref{sec: online lasso} and  \ref{sec: online debiased} are devoted to
the construction of the online lasso estimator and the online debiased method,  respectively.

\subsection{Online lasso estimator}\label{sec: online lasso}
We first consider the online lasso estimator through the gradient descent method.
Define  the score function $\bu(y_i;\bx_i,\bbeta): = \partial d(y_i;g(\bx_i^\top\bbeta))/\partial \bbeta$.
The negative gradient of the first term in \eqref{lasso_offline} $\mathcal{L}(\bbeta):=\sum_{i\in\mathcal{D}_b^\star} d(y_i;g(\bx_i^\top\bbeta))$  is
\begin{equation}\label{gradient_descent_offline}
\bU^{(b)}(\bbeta):=-\frac{\partial \mathcal{L}(\bbeta) }{\partial \bbeta}  =\sum_{i\in\mD_b^\star} \bu(y_i;\bx_i,\bbeta).
\end{equation}
 Let $U^{(j)}(\bbeta):=\sum_{i\in\mD_j}\bu(y_i;\bx_i,\bbeta)$ be the score function pertaining to the data batch $\mD_j$. Then $\bU^{(b)}(\bbeta)$ in \eqref{gradient_descent_offline} can be rewritten as a linear aggregation form:
\begin{equation*}\label{gradient_descent_online}
	\bU^{(b)}(\bbeta) =\sum_{j =1}^b \sum_{i\in\mD_j} \bu(y_i;\bx_i,\bbeta) = \sum_{j=1}^{b}U^{(j)}(\bbeta).
\end{equation*}
To pursue an online estimator upon the arrival of $\mD_b$, a key step is to update ${\bU}^{(b-1)}(\bbeta)$ to $\bU^{(b)}(\bbeta)$ without re-accessing the whole historical raw data $\mD_{b-1}^\star$.

To illustrate the idea,  we
consider a simple case with two data batches, i.e.,  $\mD_2$ arrives subsequent to $\mD_1$.
The lasso estimator based on the first batch data is denoted by $ \widehat{\bbeta}^{(1)} $, which is  the offline estimator $ \bar{\bbeta}^{(1)}$
according to (\ref{lasso_offline}).  To avoid using individual-level raw data in $\mD_1$,  we approximate ${\bU}^{(1)}(\bbeta)$ through a first-order Taylor expansion at $\widehat{\bbeta}^{(1)}$.
Note that $ \widehat{\bbeta}^{(1)} $
 satisfies $ \lVert\bU^{(1)}(\widehat{\bbeta}^{(1)})\lVert_\infty=\mathcal{O}_p(\lambda_{1}N_1)$, which implies
\begin{equation*}
\begin{split}
{\bU}^{(2)}(\bbeta) &= U^{(1)}(\bbeta) + U^{(2)}(\bbeta) \\
& = U^{(1)}(\widehat{\bbeta}^{(1)}) + \bJ^{(1)}(\widehat{\bbeta}^{(1)})(\widehat{\bbeta}^{(1)} - \bbeta) + U^{(2)}(\bbeta) + N_1\mathcal{O}_p(\|\widehat{\bbeta}^{(1)} - \bbeta\|^2),
\end{split}
\end{equation*}
where $\bJ^{(1)}(\bbeta)=-\partial U^{(1)}(\bbeta)/\partial\bbeta$ denotes the information matrix. Next, we propose  to approximate $\bU^{(2)}(\bbeta)$ by
\begin{equation}\label{eq:2_batches}
\begin{split}
\widetilde{\bU}^{(2)}(\bbeta) :=\bJ^{(1)}(\widehat{\bbeta}^{(1)})(\widehat{\bbeta}^{(1)} - \bbeta) + U^{(2)}(\bbeta).
\end{split}
\end{equation}
Apparently, calculating $\widetilde{\bU}^{(2)}(\bbeta)$ through \eqref{eq:2_batches} only requires access to the summary statistics $\{\widehat{\bbeta}^{(1)},\bJ^{(1)}(\widehat{\bbeta}^{(1)}) \}$ rather than the full dataset $\mD_1$.

Now we are ready to generalize the above approximation to an arbitrary data batch $\mD_b$. Let $\widetilde{\bJ}^{(b-1)} = \sum_{j=1}^{b-1}\bJ^{(j)}(\widehat{\bbeta}^{(j)})$ denote the aggregated information matrix. The approximation procedure in \eqref{eq:2_batches}  becomes
\begin{equation*}\label{eq:b_batches}
\begin{split}
\widetilde{\bU}^{(b)}(\bbeta)
& = \left\{\sum_{j=1}^{b-1}\bJ^{(j)}(\widehat{\bbeta}^{(j)}) \right\} (\widehat{\bbeta}^{(b-1)}-\bbeta ) + U^{(b)}(\bbeta)\\
& =\widetilde{\bJ}^{(b-1)}(\widehat{\bbeta}^{(b-1)}-\bbeta ) + U^{(b)}(\bbeta).
\end{split}
\end{equation*}
The aggregated gradient   $\widetilde{\bU}^{(b)}(\bbeta)$
depends only on $\{\widehat{\bbeta}^{(b-1)},\widetilde{\bJ}^{(b-1)}\}$. Hence,  $ \widehat{\bbeta}^{(b)} $ can be computed  through the following algorithm.
\begin{itemize}
	\item Step 1: update $\widehat{\bbeta}^{(b)}$ through
	\[
	\widehat{\bbeta}^{(b)} \leftarrow \widehat{\bbeta}^{(b)} + \frac{\eta}{2N_b} \widetilde{\bU}^{(b)}(\widehat{\bbeta}^{(b)}),
	\]
	where $\eta$ is the learning rate in the gradient descent;
	\item Step 2: for $r=1,\dots, p$, apply the soft-thresholding operator $\mathcal{S}(\widehat{\beta}_{r}^{(b)};\eta\lambda_{b})$ to the $r$-th component of $\widehat{\bbeta}^{(b)}$ obtained in Step 1:
	\[
	\mathcal{S}(\widehat{\beta}^{(b)}_{r}; \eta\lambda_b)
	=\begin{cases}
	\widehat{\beta}^{(b)}_{r} + \eta\lambda_b, \quad &\text{if} \ \widehat{\beta}^{(b)}_{r} < - \eta\lambda_b \\
	0, \quad &\text{if}\ |\widehat{\beta}^{(b)}_{r}| \leq \eta\lambda_b  \\
	\widehat{\beta}^{(b)}_{r} - \eta\lambda_b, \quad &\text{if} \ \widehat{\beta}^{(b)}_{r} > \eta\lambda_b
	\end{cases}.
	\]
\end{itemize}
Note that $\lambda_b$ is the regularization parameter that is chosen adaptively for step $b$. More details on the adaptive tuning procedure is presented in Section~\ref{sec:tuning}. 

The above two steps are carried out iteratively till convergence to obtain the online lasso estimator $ \bm{\widehat{\beta}}^{(b)} $. In the implementation, we set the
stopping criterion to be $\|{\eta}\times \widetilde{\bU}^{(b)}(\widehat{\bbeta}^{(b)})/{2N_b}\|_2 \leq 10^{-6}$.
In summary, our proposed online estimator $ \widehat{\bbeta}^{(b)} $ can be defined as
\begin{equation}\label{lasso_step_j}
	\widehat{\bbeta}^{(b)} = \underset{\bbeta\in\mathbb{R}^{p}}{\arg\min} \left[ \frac{1}{2N_b} \left\{\sum_{i\in\mathcal{D}_b} d\left(y_i; g(\bx_i^\top\bbeta)\right) + \frac{1}{2}(\bbeta - \widehat{\bbeta}^{(b-1)})^\top\widetilde{\bJ}^{(b-1)}(\bbeta - \widehat{\bbeta}^{(b-1)})    \right\} + \lambda_{b}\lVert \bbeta \lVert_1  \right].
\end{equation}
In contrast to  the standard offline lasso estimator in \eqref{lasso_offline}, our proposed online estimator in \eqref{lasso_step_j} depends on the data only through the summary statistics $\{\widehat{\bbeta}^{(b-1)},\widetilde{\bJ}^{(b-1)} \}$ instead of $\mD_{b-1}^\star$.

\subsection{Online debiased lasso}  \label{sec: online debiased}
We now proceed to study the online statistical inference and construct confidence intervals for the $r$-th component of the regression parameter vector, $r=1,\ldots,p$. We first define a low-dimensional projection
\begin{equation}\label{projection_step_j}
	\bm{\widehat{\gamma}}_r^{(b)}= \underset{\bm{\gamma}\in\mathbb{R}^{(p-1)}}{\arg\min}
	\left\{\frac{1}{2N_b}\left(\widetilde{\bJ}^{(b)}_{r,r} - 2 \widetilde{\bJ}^{(b)}_{r,-r}\bgamma +\bgamma^\top\widetilde{\bJ}_{-r,-r}^{(b)}\bgamma \right) +\lambda_b\|\bm{\gamma} \|_1
	\right\}.
\end{equation}
where $ N_b $ and $ \lambda_{b} $ are the same as in \eqref{lasso_step_j}.
Letting $ \widehat{W}{^{(j)}} \in \mathbb{R}^{n_j\times n_j}$ be the diagonal matrix with diagonal elements $ \sqrt{g'(\bm{X}^{(j)}\bm{\widehat{\beta}}^{(j)} )} $ and $ \widehat{\bX}^{{(j)}} $ be the weighted design matrix $ \widehat{W}{^{(j)}}\bm{X}^{(j)}$, then  \eqref{projection_step_j} can be recast into
\begin{equation}\label{projection_step_j_version2}
	\bm{\widehat{\gamma}}_r^{(b)}= \underset{\bm{\gamma}\in\mathbb{R}^{(p-1)}}{\arg\min}
\left\{\frac{1}{2N_b}\sum_{j=1}^{b}\left\lVert \widehat{\bX}^{{(j)}}_{\cdot r}  - \widehat{\bX}^{{(j)}}_{\cdot -r}\bgamma  \right\lVert_2^2 +\lambda_b\|\bm{\gamma} \|_1
\right\}.
\end{equation}
With this view, $ \bm{\widehat{\gamma}}_r^{(b)}  $  can be computed in a similar fashion as
 the online lasso estimator.
 Specifically, \eqref{projection_step_j_version2} has the same form as \eqref{lasso_offline} if we choose the deviance  function $ d(\cdot;\cdot) $ as the square error. Then, it is straightforward to propose the online algorithm to find $ \bm{\widehat{\gamma}}_r^{(b)} $ by following the procedure in Section \ref{sec: online lasso}. At this time, the summary statistics is
$ (\widetilde{\bJ}^{(b)}_{r,-r},\widetilde{\bJ}_{-r,-r}^{(b)}), $ which has been already stored in $ \widetilde{\bJ}^{(b)}. $  As a result, we find a solution of \eqref{projection_step_j}.
 After obtaining $\widehat{\bgamma}_r^{(b)}$ in \eqref{projection_step_j}, we compute
$$
\widehat{\bz}_r^{(b)} : = \bx_r^{(b)} - \bX_{-r}^{(b)} \widehat{\bgamma}_r^{(b)}, \ \widehat{\tau}_r^{(b)}:=\widetilde{\bJ}_{r,r}^{(b)} - \widetilde{\bJ}_{r,-r}^{(b)} \widehat{\bgamma}_r^{(b)}.
$$
Denote $ \widetilde{\bgamma}^{(j)}_r  = (\widehat{\bgamma}^{(j)}_{r,1}, \ldots, -1, \ldots, \widehat{\bgamma}^{(j)}_{r,p} )^\top \in \mathbb{R}^p$, whose $ r$-th element is $ -1 $.
Then,
upon the arrival of the batch data $\mD_b$, we define our proposed online debiased lasso estimator as
\begin{align}\label{online_debiased_algorithm}
{\widehat{\beta}}^{(b)}_{\text{on}, r} &= \widehat{{\beta}}^{(b)}_{r} + \frac{1}{\widehat{\tau}_r^{(b)}}    \left[\sum_{j=1}^{b}\{\bm{\widehat{z}}^{(j)}_{r}\}^\top\left\{\by^{(j)} - g(\bm{X}^{(j)}\bm{\widehat{\beta}}^{(j)} )\right\} + \sum_{j=1}^{b}\{\widetilde{\bgamma}^{(j)}_r\}^\top\bJ^{(j)}(\widehat{\bbeta}^{(j)})\{{\widehat{\bbeta}}^{(b)} - {\widehat{\bbeta}}^{(j)} \}\right] \nonumber \\
&\equiv  \widehat{{\beta}}^{(b)}_{r} + \text{debiasing term} + \text{online error correction term}.
\end{align}
It is worth noting that the debiased lasso estimator
involves
 the initial lasso estimator defined in \eqref{lasso_step_j}, as well as two additional terms: {\it a debiasing term} and {\it an online error correction term}. \cite{vandegeer2014} studied the offline version of debiased lasso for GLM. Our debiased term could be viewed as an online generalization of the offline counterpart in \cite{vandegeer2014}.
 However, they are rather different in the sense that,   the debiased term alone  will not suffice to ensure the asymptotic normality. As  $ \widetilde{\bJ}^{(b)} $
 is used to approximate the information matrix $ \bJ(\bbeta^0)$, the approximation error accumulates even though each $ {\widehat{\bbeta}}^{(j)}, j = 1,\ldots, b, $ is  consistent for $ \bbeta^{0} $. The correction term in (\ref{online_debiased_algorithm}) is used to eliminate the approximation error arising from the online update. 

Meanwhile, the proposed debiased lasso estimator with the
online error correction term aligns with the online learning framework, as
  \eqref{online_debiased_algorithm}  only requires the following summary statistics rather than the entire dataset $ \mD_{b}^\star $:
\begin{align}
		s_{1}^{(b)} &= \sum_{j=1}^{b}\{\bm{\widehat{z}}^{(j)}_{r}\}^\top\left(\bm{y}^{(j)} - g(\bX^{(j)}\widehat{\bbeta}^{(j)})\right), \nonumber\\
		s_{2}^{(b)} &=  \sum_{j=1}^{b} \{\widetilde{\bgamma}^{(j)}_r\}^\top\bJ^{(j)}(\widehat{\bbeta}^{(j)}) {\widehat{\bbeta}}^{(j)},\ \
		\bS^{(b)}  =\sum_{j=1}^{b} \{\widetilde{\bgamma}^{(j)}_r\}^\top\bJ^{(j)}(\widehat{\bbeta}^{(j)}), \label{other statistics}
\end{align}
which keep the same size when  new data arrive, and can be easily updated.

The asymptotic normality of the online debiased lasso and the consistency of two lasso-typed estimators in \eqref{lasso_step_j} and \eqref{projection_step_j}
are established  in Section \ref{sec: asymptotic}. 
For variance estimation,  let
\begin{equation}
v^{(b)} =\sum_{j=1}^{b} \{\bm{\widehat{z}}^{(j)}_{r}\}^\top
\left\{\bm{y}^{(j)} - g(\bX^{(j)}\widehat{\bbeta}^{(j)})\right\} \left\{\bm{y}^{(j)} - g(\bX^{(j)}\widehat{\bbeta}^{(j)})\right\}^\top \bm{\widehat{z}}^{(j)}_{r}. \label{other stat: var}
\end{equation}
The estimated standard error
$\widehat{\sigma}^{(b)}_{r} := {\sqrt{v^{(b)}}}/{\widehat{\tau}_r^{(b)}}$ can also be updated online accordingly.

\subsection{Adaptive tuning}\label{sec:tuning}
In an offline setting, the regularization parameter $\lambda$ is typically determined by cross-validation where the entire dataset is split into training and test sets multiple times. However, since the full dataset is not accessible in an online setting, such procedure is not feasible. To align with the nature of streaming datasets, we use the ``rolling-original-recalibration" procedure with mean squared prediction error (MSE) being the cross-validation criterion~\citep{han2021online}. Specifically, at time point $b$, cumulative dataset up to time point $b-1$ serves as the training set while the new data batch $\mD_b$ is the test set. It is worth noting that instead of re-accessing raw data $\{\mD_1,\dots,\mD_{b-1}\}$, we plug in $\widehat{\bbeta}^{(b-1)}$ to evaluate the prediction error for a sequence of $\lambda$ in a candidate set $T_{\lambda}$:
\begin{equation}\label{eq:MSE}
PE_b(\lambda)=n_b^{-1}\|\by^{(b)} - g(\bX^{(b)}\widehat{\bbeta}^{(b-1)}(\lambda)) \|^2_2,\ \lambda\in T_{\lambda},
\end{equation}
and choose $\lambda$ such that $\lambda_b:=\underset{\lambda\in T_{\lambda}}{\arg\min}\ PE_b(\lambda)$. The initial $\lambda_1$ is selected by the classical offline cross-validation.

\subsection{Algorithm}\label{sec:algorithm}
\begin{figure}[h]
	\centering
	\includegraphics[width=\linewidth]{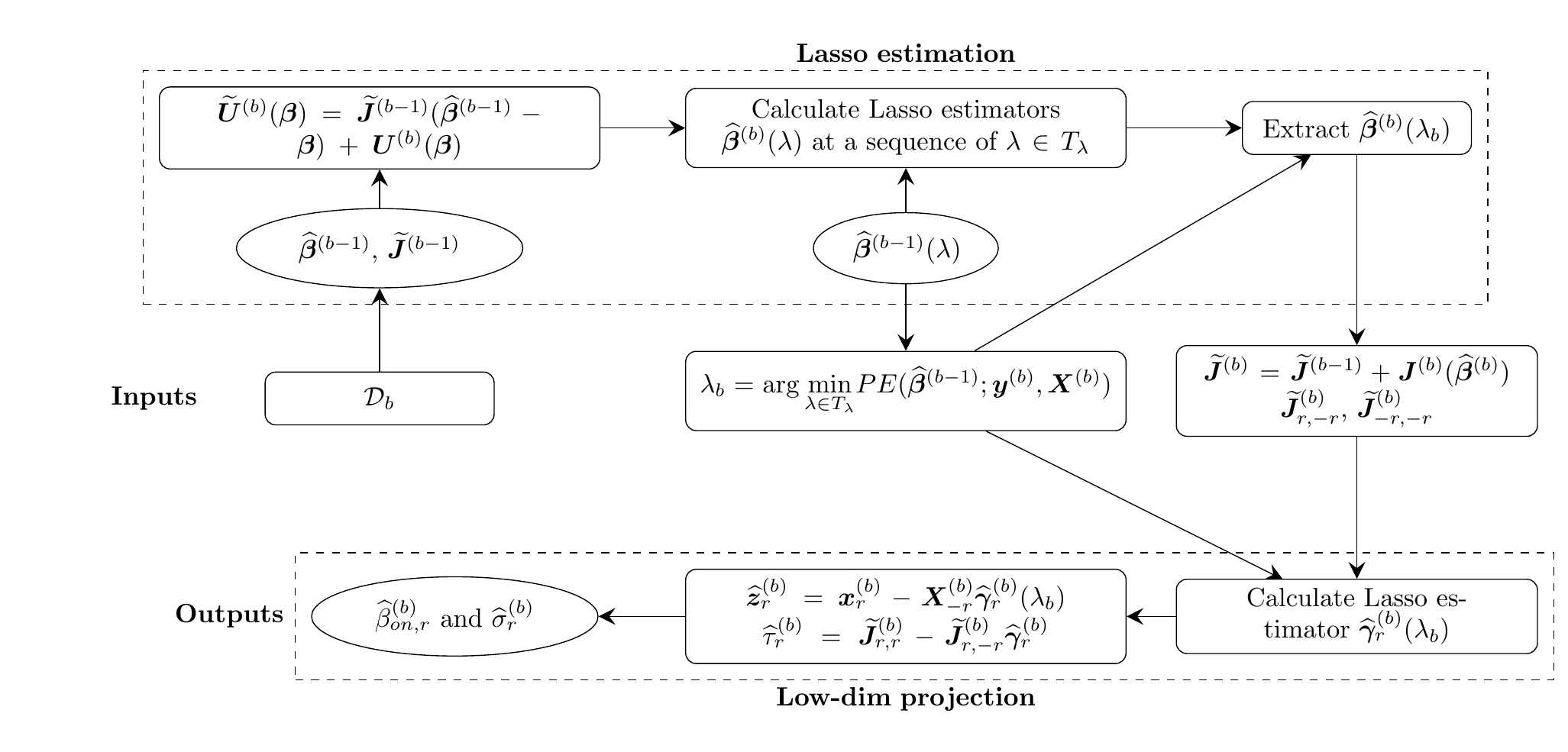}
	\caption{\label{fig:algorithm}Flowchart of the online debiasing algorithm. When a new data batch $\mD_b$ arrives, it is sent to the  \textit{lasso estimation} block for updating $\widehat{\bbeta}^{(b-1)}$ to $\widehat{\bbeta}^{(b)}$. At the same time, it is also viewed as test set for adaptively choosing tuning parameter $\lambda_b$. In the \textit{low-dim projection}  block, we extract sub-matrices from the updated information matrix $\widetilde{\bJ}^{(b)}$ to compute $\widehat{\bgamma}_r^{(b)}(\lambda_b)$ and the corresponding low-dimensional projection $\widehat{\bz}_r^{(b)}$ and $\widehat{\tau}_r^{(b)}$. Outputs $\widehat{\beta}_r^{(b)}(\lambda_b)$, $\widehat{\bz}_r^{(b)}$ and $\widehat{\tau}_r^{(b)}$ are further used to compute the debiased lasso estimator $\widehat{\beta}_{\text{on},r}^{(b)}$ and its estimated standard error $\widehat{\sigma}_{r}^{(b)}$.}
\end{figure}

We summarize the results in Sections \ref{sec: online lasso}-\ref{sec:tuning} in Figure~\ref{fig:algorithm} and Algorithm \ref{alg: online_debiased_lasso}. It consists of two main blocks: one is~\textit{online lasso estimation} and the other is \textit{online low-dimensional projection}. Outputs from both blocks are used to compute the online debiased lasso estimator as well as the construction of confidence intervals in real-time. In particular, when a new data batch $\mathcal{D}_b$ arrives,
it is first sent to the \textit{online lasso estimation} block, where the summary statistics $\left\{\widehat{\bbeta}^{(b-1)},\widetilde{\bJ}^{(b-1)} \right\}$ are used to compute $\widetilde{\bU}^{(b)}$. Then we use gradient descent to update the lasso estimator $\widehat{\bbeta}^{(b-1)}$ to $\widehat{\bbeta}^{(b)}$
at a sequence of tuning parameters values without retrieving the whole dataset. At the same time, regarding the cumulative dataset that produces the old lasso estimate $\widehat{\bbeta}^{(b-1)}$ as training set and the newly arrived $\mathcal{D}_b$ as testing set, we can choose the tuning parameter $\lambda_b$ that gives the smallest prediction error. Now, the selected $\lambda_b$ and sub-matrices of $\widetilde{\bm{J}}^{(b)}$ are passed to the \textit{low-dimensional projection} block for the calculation of  $\widehat{\bgamma}_r^{(b)}(\lambda_b)$. The resulting projection $\widehat{\bz}_r^{(b)}$ and $\widehat{\tau}_r^{(b)}$ from the \textit{low-dimensional projection} block together with the lasso estimator $\bm{\widehat{\beta}}^{(b)}$ will be used to compute the debiased lasso estimator $\widehat{\beta}^{(b)}_{\text{on},r}$ and its estimated standard error $\widehat{\sigma}_r^{(b)}$.

\begin{algorithm}
	\caption{Online Debiased Algorithm in GLM}
	\centering
	\label{alg: online_debiased_lasso}
	\begin{algorithmic}[5]
		\Require
		Learning rate $ \eta $, a set of candidate values $T_{\lambda}$ for penalty parameter;
		\Ensure
		Debiased lasso estimator $\widehat{\beta}_{\text{on}, r}^{(b)}$;
		\For{$b = 1, 2,  \ldots $}
		\State Receive the streaming dataset $ \mD_{b} $;
		\State For a sequence of $\lambda\in T_{\lambda}$, update lasso estimator $ \widehat{\bbeta}^{(b)}(\lambda) $ defined in \eqref{lasso_step_j} via gradient descent (GD);
		\State Determine $\lambda_b$ based on minimum MSE, $\lambda_b:=\underset{\lambda\in T_{\lambda}}{\arg\min} \ PE_b(\lambda)$ defined in~\eqref{eq:MSE};
		\State Update and store the summary statistics $ \{\widehat{\bbeta}^{(b)}, \widetilde{\bJ}^{(b)}\} $;
		\State Given $\lambda_b$, find the low-dimensional projection $ \widehat{\bgamma}_r^{(b)} $ defined in \eqref{projection_step_j} via GD;
		\State Update and store the summary statistics $ s_1^{(b)}, s_2^{(b)} $ and $ \bS^{(b)} $ by \eqref{other statistics};
		\State Compute $\widehat{\beta}_{\text{on}, r}^{(b)}$ by \eqref{online_debiased_algorithm} and $ \widehat{\sigma}^{(b)}_{r} $;
		\EndFor \State {\textbf{end for}}
		
		\Return $ \widehat{\beta}_{\text{on}, r}^{(b)} $ and its estimated standard error $ \widehat{\sigma}^{(b)}_{r} $.
	\end{algorithmic}
\end{algorithm}

\begin{remark}
	Algorithm \ref{alg: online_debiased_lasso} computes  the estimators in the chronological order, namely \\
	$ (\widehat{\bbeta}^{(1)}, \widehat{\bgamma}_r^{(1)}, \widehat{\beta}_{\text{on}, r}^{(1)}, \ldots, \widehat{\bbeta}^{(b)}, \widehat{\bgamma}_r^{(b)}, \widehat{\beta}_{\text{on}, r}^{(b)}) $. Since  $  \widetilde{\bJ}^{(b-1)} $ involves $\widehat{\bbeta}^{(1)}, \ldots, \widehat{\bbeta}^{(b-1)} $, there are dependency among the estimators. For instance, both
	$ \widehat{\bgamma}_r^{(b-1)} $  and $ \widehat{\bbeta}^{(b)}$   depend on the previous lasso estimation $ (\widehat{\bbeta}^{(1)},\ldots, \widehat{\bbeta}^{(b-1)}). $
\end{remark}

\begin{remark}
		When $p$ is 
	large, 
	it may be challenging to implement Algorithm \ref{alg: online_debiased_lasso} 
	since the space complexity
	to store the aggregated information matrix $\widetilde{\bJ}^{(b)}$ is $\mathcal{O}(p^2).$
	To reduce memory usage, we can compute the eigenvalue decomposition (EVD) of $\widetilde{\bJ}^{(b)}=Q_{b}\Lambda_{b} Q_{b}^\top$, where $Q_{b}$ is the $p\times N_{b}$ columns orthogonal matrix of the eigenvectors, $\Lambda_{b}$ is the $N_{b}\times N_{b}$ diagonal matrix whose diagonal elements are the eigenvalues of $\widetilde{\bJ}^{(b)}$. We only need to store $Q_{b}$ and $\Lambda_{b}.$
		Since
		$ r_{b}=  $ rank$ (\Lambda_{b}) \leq \min\{N_{b}, p \},$ 
		we can use an incremental EVD approach~\citep{onlinePCA2018} to update $Q_{b}$ and $\Lambda_{b}$. Then the space complexity reduces to $\mathcal{O}(r_{b}p)$. The space complexity can be further reduced by setting a threshold. For example, select the principal components which explain most of the variations in the predictors. 
		However, incremental EVD 
		could increase the computational cost since it requires additional $ \mathcal{O}(r_{b}^2p) $ computational complexity. Indeed, there is a trade-off between the space complexity and computational complexity. How to balance this trade-off is an important computational issue
		and deserves careful analysis,  but is beyond the scope of this study.
\end{remark}

\section{Asymptotic properties}\label{sec: asymptotic}
In this section, we  state  our main theoretical results: the consistency of lasso estimators $ \widehat{\bbeta}^{(b)} $ and  $ \bm{\widehat{\gamma}}_r^{(b)} $ defined in \eqref{lasso_step_j} and \eqref{projection_step_j} respectively, as well as the asymptotic normality of the online debiased estimator $ {\widehat{\beta}}^{(b)}_{\text{on}, r}. $

Before stating the main  results,  some notations are needed. Recall that $ \bbeta^0 $ is the true coefficient. Consider a random design matrix $ \bX $ with  i.i.d rows. Let $ \bJ $ denote the population version of the information matrix in GLM, i,e., $ \bJ = \mathbb{E} [\bJ^{(1)}(\bbeta^0)]/N_1, $ and let $ \bTheta  = \bJ^{-1}$ be its inverse. Then, the ground truth of $ \bm{\widehat{\gamma}}_r^{(b)}$ given in \eqref{projection_step_j} is defined as
\begin{eqnarray*}
	\bm{{\gamma}}_r^{0}= \underset{\bm{\gamma}\in\mathbb{R}^{(p-1)}}{\arg\min}\ \mathbb{E}
	\left({\bJ}_{r,r} - 2 {\bJ}_{r,-r}\bgamma +\bgamma^\top {\bJ}_{-r,-r}\bgamma
	\right).
\end{eqnarray*}
Recall that $S_0:=\{r:\beta^0_{r}\neq 0\}$  and $s_0 = \lvert S_0 \lvert$. In addition, $ S_r =  \{k \neq r: \bTheta_{k,r} \neq 0  \} $ and $ s_r = \lvert S_r \lvert$ for $r = 1, \ldots, p. $
 The following assumptions are needed to establish the consistency of the lasso estimators $ \widehat{\bbeta}^{(b)} $ and  $ \widehat{\bgamma}_r $ defined in \eqref{lasso_step_j} and \eqref{projection_step_j} respectively.
\begin{assumption}\label{assumption_1}
Suppose that
	
	(A1) The pairs of random variables $ \{y_i, \bx_i \}_{i \in \mathcal{D}_b^\star} $ are i.i.d..  The covariates are bounded by some finite constant $ K > 0$, i.e.,  $ \sup_{i \in \mathcal{D}_b^\star} \lVert \bx_i \lVert_\infty \leq K.$

	(A2) $ \sup_{i \in \mathcal{D}_b^\star} \lvert \bx_i\bbeta^0 \lvert = \mathcal{O}(1)$ and $\sup_{i \in \mathcal{D}_b^\star} \lvert (\bx_{i})_{-r}\bm{{\gamma}}_r^{0} \lvert = \mathcal{O}(1),$ where $ (\bx_{i})_{-r} $ is the sub-vector of $ \bx_{i} $ with $ r$-th element deleted.
		
		
	(A3) The derivative $ g'(\cdot) $ exists. For some $ \delta$-neighborhood $ ( \delta > 0 )$, $ g'(\cdot) $ is  locally Lipschitz with constant $ L $, that is,
	$$\sup_{i \in \mathcal{D}_b^\star}\ \sup_{\bbeta,\bbeta' \in \{\bbeta: \lVert \bbeta - \bbeta^0 \lVert_1 \leq \delta\}} \frac{|g'(\bx_{i}\bbeta) - g'(\bx_{i}\bbeta')|}{|\bx_{i}\bbeta-\bx_{i}\bbeta'|} \leq L. $$	
	In addition, for all $ \lVert \bbeta - \bbeta^0 \lVert \leq \delta$, $ \sup_{i \in \mathcal{D}_b^\star} \lvert {g'(\bx_i \bbeta)} \lvert = \mathcal{O}(1)$ and $ \sup_{i \in \mathcal{D}_b^\star} \lvert {1/g'(\bx_i \bbeta)} \lvert = \mathcal{O}(1). $

	{ (A4) The smallest eigenvalue of $ \bJ $ is bounded away from zero.}
	
\end{assumption}
The above assumptions are regular conditions  in studying the properties of lasso \citep{van2008high,vandegeer2014}.
(A1) assumes the streaming data is homogeneous and the covariates are sampled from some bounded distribution.  (A2) is imposed to avoid some extreme cases.
(A3) requires the local Lipschitz property of the derivative of the mean function around the  truth value. It can be easily verified that the popular logistic regression, a special case of GLMs, satisfies this condition.
(A4) is needed to ensure that the compatibility condition holds.
\vspace{0.1in}

\begin{theorem}\label{thm: lasso consistency}
	Assume Assumption \ref{assumption_1} holds.
	Suppose that the first batch size $ n_1 \geq c s_0\log p $ for some constant $ c $ and $ b = o(\log N_b) $, and the tuning parameter $ \lambda_j = C\sqrt{\log p/ N_j}, j = 1,\ldots b$ for some constant $ C $.
	Then, for any $ j = 1, \ldots, b $, with probability at least $  1 - p^{-2} $,
	the proposed online estimator	
	in \eqref{lasso_step_j} satisfies
	\begin{equation}\label{lasso_consistency}
	{
		\lVert {\widehat{\bbeta}}^{(j)} - {{\bbeta}}^{0} \lVert_1 \leq c_{1}^{(j)} s_0 \lambda_j, \ \lVert \bX^{(j)}_\star({\widehat{\bbeta}}^{(j)} - {{\bbeta}}^{0}) \lVert^2_2 \leq c_{2}^{(j)} s_0N_j \lambda_j^2.
		}
	\end{equation}
\end{theorem}

\vspace{0.2cm}

\begin{remark}
	Theorem \ref{thm: lasso consistency} provides upper bounds of the estimation error and the prediction error of the online lasso estimator $ \widehat{\bbeta}^{(b)} $. The constants $ c_1^{(j)} $ and $ c_2^{(j)}, j = 1, \ldots, b, $ possibly depend on the batch step $ b $. Recall that the lasso estimator $ \widehat{\bbeta}^{(b)} $
	depends on $ \widehat{\bbeta}^{(b-1)} $. So the estimation error in the previous step will be carried onto  the updated estimators. As a result, it is
	inevitable that some constants in the oracle inequality depend on
	 $ b $; nonetheless, they are well under control  as long as $ b = o(\log N_b). $
\end{remark}
The next corollary shows the consistency of the proposed online lasso estimator in \eqref{lasso_step_j}.
\begin{corollary}\label{cor: consistency}
	Assume those conditions in Theorem \ref{thm: lasso consistency} hold. If there exists an $ \epsilon > 0 $ such that
	$ s_0^2\log(p) N_b^{-{1}+\epsilon} = o(1), $
then the lasso estimator in \eqref{lasso_step_j} satisfies
	\begin{equation*}
			\lVert {\widehat{\bbeta}}^{(b)} - {{\bbeta}}^{0} \lVert_1 \rightarrow_p 0 \ \text{as} \ N_b\to \infty,
	\end{equation*}
	where  $\rightarrow_p$ means convergence in probability.
	\end{corollary}
	
Similarly, we present the oracle inequality for $ \widehat{\bgamma}_{r}^{(b)} $ in the next theorem.

\begin{theorem}\label{thm: lasso consistency 2}
	Assume Assumption \ref{assumption_1} holds and  the cumulative batch size satisfies $ N_j \geq cc^{(j)}_1 s_0^2s_r^2\log p, j = 1, \ldots, b$ for some constant $ c $.
	Then, for any $ j = 1, \ldots, b $, with probability at least $  1 - p^{-2} $, low-dimensional projection defined in \eqref{projection_step_j} with $ \lambda_j =
	c\sqrt{{\log p}/{N_j}} $ satisfies
	\begin{equation}\label{gamma_consistency}
		\lVert {\widehat{\bgamma}_r}^{(j)} - {{\bgamma}}_{r}^{0} \lVert_1 \leq c_3 s_r \lambda_j.
	\end{equation}
\end{theorem}

Combining the results in Theorem \ref{thm: lasso consistency} and Theorem \ref{thm: lasso consistency 2}, we are ready to establish the asymptotic normality of the ODL estimator. 
\begin{theorem}\label{thm: lasso asymptotic normality}
Assume those conditions in Theorem \ref{thm: lasso consistency} and Theorem \ref{thm: lasso consistency 2} hold. If there exists an $ \epsilon > 0 $ such that
\begin{equation*}
 s_0^2\log(p)\log(N_b) N_b^{-\frac{1}{2}+\epsilon} = o(1), \ s_0s_r\log(p)N_b^{-\frac{1}{2}+\epsilon} = o(1),
\end{equation*}
then for sufficiently large $ N_b $,
\begin{eqnarray*}
	&&\frac{\widehat{\tau}_r^{(b)}}{\sqrt{N_b}}({\widehat{\beta}}^{(b)}_{\text{on}, r} - \beta_{r}^{0}) = W_r + V_r,\\
	&& W_r = \frac{1}{\sqrt{N_b}}\sum_{j=1}^{b} \{{\widehat{\bz}}^{(j)}_{r}\}^\top\left(\by^{(j)} - g(\bX^{(j)}{{\bbeta}^0} )\right), V_r = o_p(1).
\end{eqnarray*}
\end{theorem}
According to Theorem 3, the asymptotic expression of $ {\widehat{\tau}_r^{(b)}}({\widehat{\beta}}^{(b)}_{\text{on}, r} - \beta_{r}^{0})/{\sqrt{N_b}} $ is a sum of
$W_r$ and $V_r$, where  $ W_r $ converges in distribution to a normal distribution by the martingale central limit theorem and $ V_r $ diminishes as $ N_b $ goes to infinity.
\begin{remark}
	 Theorem \ref{thm: lasso asymptotic normality} implies  that the total data size $ N_b $ could be as small as the logarithm of the dimensionality $ p $, which is a common condition for offline debiased lasso in the literature \citep{Zhang_delasso_2014, vandegeer2014}. However, the setting of the offline debiased lasso is quite different from the online framework considered in our paper. Due to the lack of access to the whole dataset, it is increasingly difficult to derive the asymptotic property of online debiased lasso. One major difficulty is that, we approximate the information matrix by $  \widetilde{\bJ}^{(b)} $, which involves estimators obtained in the previous steps. As a result, there is dependence among $\widehat{\bbeta}^{(1)}, \ldots, \widehat{\bbeta}^{(b-1)}, \widehat{\bbeta}^{(b)}. $
	 Another difficulty is to  deal with the approximation error accumulates in the online updating, especially under high-dimensional settings.  In contrast,
	 the classical offline lasso does not have these two problems.
Even for the online debiased lasso in linear model \citep{han2021online},
the above two problems can be bypassed by making use of the special structure of the least squares in linear model.
\end{remark}

\section{Simulation studies}\label{sec:sim}
\subsection{Setup}
In this section, we conduct simulation studies to examine the finite-sample performance of the proposed
ODL procedure in high-dimensional GLMs.
We randomly generate a total of $N_b$ samples arriving in a sequence of $b$ data batches, denoted by $\{\mathcal{D}_1,\dots,\mathcal{D}_b\}$, from the logistics regression model. Specifically,
\[
\mathbb{E}(y^{(j)}_i\mid \bx_i^{(j)}) = \frac{\exp(\{\bx_i^{(j)}\}^\top\bbeta^0)}{1 + \exp(\{\bx_i^{(j)}\}^\top\bbeta^0)}, \ i=1,\dots,n_j; \ j = 1,\dots,b,
\]
where $\bx_i^{(j)}\sim \mathcal{N}(\bm{0},\bSigma)$ 
and $\bm{\beta}^0\in\mathbb{R}^p$ is a $p$-dimensional sparse parameter vector.
Recall that $s_0$ is the number of non-zero components of $\bbeta^0$.
We set half of the nonzero coefficients to be $1$ (relatively strong signals),  and another half to be $0.01$ (weak signals).
We consider the following settings: (i) $N_b=120$, $b=12$, $n_j=10$ for  $j=1, \ldots, 12$, $p=100$ and $s_0=6$; (ii) $N_b=624$, $b=12$, $n_j=52$ for $j=1,\ldots,12$, $p=600$ and $s_0=10$. Under each setting,  two types of  $\bSigma$ are considered:
(a) $\bSigma=\bI_p$; (b) $\bSigma = \{0.5^{|i-j|} \}_{i,j=1,\dots,p}$. We set the step size in gradient descent  $\eta=0.005$.

The objective is to conduct both estimation and inference along the arrival of a sequence of data batches. The evaluation criteria include:  averaged absolute bias in estimating $\bbeta^0$ (A.bias);  averaged estimated standard error (ASE);  empirical standard error (ESE); coverage probability (CP)
of the 95\% confidence intervals;  averaged length of the 95\% confidence interval (ACL). These metrics will be evaluated separately for three groups: (i) $\beta^0_{r}=0$, (ii) $\beta^0_{r} = 0.01$ and (iii) $\beta^0_{r}=1$. Comparison is made among (i) the maximum likelihood estimator obtained by fitting the conventional GLM at the terminal point $b$ where $N_b>p$, (ii) the offline debiased $\ell_1$-penalized estimator at the terminal point $b$ which is also the benchmark method, and (iii) our proposed ODL 
estimator at several intermediate points from $j=1,\dots,b$. We include two offline methods in comparison using R packages~\texttt{hdi}~\citep{hdi2015} and~\texttt{glm}, respectively. The results are reported in Tables~\ref{tab:p_100_ind}-\ref{tab:p_600_ar1}.

\subsection{Bias and coverage probability}
It can be seen from Tables~\ref{tab:p_100_ind}-\ref{tab:p_600_ar1} that the estimation bias of the ODL 
estimator decreases rapidly as the number of data batches $b$ increasing from $2$ to $12$. Both the estimated standard errors and averaged length of 95\% confidence intervals exhibit similar decreasing trend over time, and almost coincide with those by the offline benchmark method at the terminal points. Moreover,
ODL enjoys great computation efficiency and it is almost 16 times faster than its offline counterpart.

It is worth noting that even though at the terminal point where the cumulative sample size $N_b$ is slightly larger than $p$, we can still fit the conventional GLM to obtain the MLE. It fails to provide reliable coverage probabilities due to severely large biases and estimated standard errors. In particular, the estimation bias of MLE becomes 50 times that of offline or online debiased lasso when $p=100$ as shown in Tables~\ref{tab:p_100_ind} and~\ref{tab:p_100_ar1}, and it further increases to 80 times the debiased estimators when $p=600$. Furthermore, as clearly indicated by the large empirical standard errors, MLE under this setting suffers from severe instability. Such an invalid estimation and inference results by MLE further demonstrates the advantage of our proposed online debiased method under the high-dimensional sparse logistic regression setting with streaming datasets.



\begin{table}
	\caption{\label{tab:p_100_ind} $N_b=120$, $b=12$, $p=100$, $s_0=6$, $\bSigma=\bI_p$. Performance on statistical inference. ``MLE" is the offline estimator obtained by fitting the traditional GLM, ``offline" corresponds to the offline debiased $\ell_1$-norm penalized estimator, and ``ODL" represents our proposed online debiased lasso estimator. Tuning parameter $\lambda$ is chosen from $T_\lambda= \{10^{-4}, 10^{-3}, 0.01, 0.05 \}$. Simulation results are summarized over 200 replications. In the table, we report the $\lambda$ selected with highest frequency among 200 replications.}
	\centering
	\begin{tabular}{l l| r| r | rrrrrr}
		&$\beta_{0,r}$ &MLE  &offline &\multicolumn{6}{c}{ ODL} \\
		data batch index &&&&2  &4 &6 &8  &10 &12\\
		$\lambda$ &&&&0.05 &0.05 &0.05 &0.05 &$10^{-4}$ &$10^{-4}$\\
		\hline
		\hline
		\multirow{3}{*}{A.bias} &0 
		&2.290 &0.035
		&0.080 &0.056 &0.049 &0.038 &0.036 &0.033\\
		&0.01 
		&1.150 &0.038
		&0.131 &0.036 &0.035 &0.033 &0.020 &0.034\\
		&1 
		&17.87 &0.057
		&0.176 &0.125 &0.125 &0.125 &0.120 &0.109\\
		\hline
		\multirow{3}{*}{ASE}  &0 
		&$1.78\times10^6$ &0.595
		&1.336 &0.949 &0.777 &0.679 &0.609 &0.559\\
		&0.01  
		&$1.79\times10^6$ &0.591
		&1.323 &0.939 &0.771 &0.674 &0.604 &0.555\\
		&1    
		&$1.83\times10^6$ &0.598
		&1.341 &0.950 &0.778 &0.680 &0.609 &0.561\\
		\hline
		\multirow{3}{*}{ESE}  &0 
		&34.90 &0.591
		&1.367 &0.962 &0.786 &0.681 &0.611 &0.559\\
		&0.01 
		&36.25 &0.581
		&1.341 &0.996 &0.792 &0.686 &0.597 &0.553\\
		&1  
		&33.49 &0.585
		&1.327 &0.903 &0.755 &0.664 &0.593 &0.548\\
		\hline
		\multirow{3}{*}{CP}   &0 
		&1.000 &0.953
		&0.951 &0.949 &0.948 &0.949 &0.950 &0.950\\
		&0.01  
		&1.000 &0.960
		&0.943 &0.938 &0.943 &0.947 &0.960 &0.953\\
		&1  
		&1.000 &0.965
		&0.955 &0.957 &0.962 &0.953 &0.942 &0.953\\
		\hline
		\multirow{3}{*}{ACL}   &0 
		&$6.96\times10^6$ &2.487
		&5.239 &3.721 &3.047 &2.663 &2.385 &2.190\\
		&0.01  
		&$7.02\times10^6$ &2.478
		&5.187 &3.681 &3.022 &2.641 &2.366 &2.176\\
		&1  
		&$7.18\times10^6$ &2.479
		&5.258 &3.725 &3.051 &2.666 &2.387 &2.197\\
		\hline
		C.Time (s)  &&0.03 &29.97  &\multicolumn{6}{c}{2.026}\\
	\end{tabular}
\end{table}

\begin{table}
	\caption{\label{tab:p_100_ar1} $N_b=120$, $b=12$, $p=100$, $s_0=6$, $\bSigma= \{0.5^{|i-j|} \}_{i,j=1,\dots,p}$. Performance on statistical inference. ``MLE" is the offline estimator obtained by fitting the traditional GLM, ``offline" corresponds to the offline debiased $\ell_1$-norm penalized estimator, and ``ODL" represents our proposed online debiased lasso estimator. Tuning parameter $\lambda$ is chosen from $T_\lambda= \{10^{-4}, 10^{-3}, 0.01, 0.05 \}$. Simulation results are summarized over 200 replications. In the table, we report the $\lambda$ selected with highest frequency among 200 replications.}
	\centering
	\begin{tabular}{l l| r| r | rrrrrr}
		&$\beta_{0,r}$ &MLE  &offline &\multicolumn{6}{c}{ODL} \\
		data batch index &&&&2  &4 &6 &8  &10 &12\\
		$\lambda$ &&&&0.05 &0.05 &0.05 &0.05 &$10^{-4}$ &$10^{-4}$\\
		\hline
		\hline
		\multirow{3}{*}{A.bias} &0 
		&2.467 &0.053
		&0.119 &0.093 &0.084 &0.078 &0.073 &0.068\\
		&0.01 
		&1.546  &0.048
		&0.138 &0.070 &0.093 &0.068 &0.065 &0.064\\
		&1 
		&17.54  &0.026
		&0.122 &0.122 &0.113 &0.087 &0.095 &0.084\\
		\hline
		\multirow{3}{*}{ASE}  &0 
		&$2.45\times10^6$ &0.674
		&1.341 &0.955 &0.781 &0.680 &0.614 &0.567\\
		&0.01  
		&$2.43\times10^6$ &0.673
		&1.332 &0.948 &0.780 &0.679 &0.612 &0.565\\
		&1    
		&$2.36\times10^6$ &0.667
		&1.344 &0.962 &0.783 &0.682 &0.616 &0.567\\
		\hline
		\multirow{3}{*}{ESE}  &0 
		&45.93 &0.671
		&1.358 &0.953 &0.781 &0.679 &0.608 &0.556\\
		&0.01 
		&45.47 &0.667
		&1.387 &0.935 &0.765 &0.644 &0.595 &0.557\\
		&1  
		&45.87 &0.657
		&1.334 &0.940 &0.774 &0.659 &0.592 &0.544\\
		\hline
		\multirow{3}{*}{CP}   &0 
		&1.000 &0.951
		&0.951 &0.951 &0.948 &0.947 &0.947 &0.949\\
		&0.01  
		&1.000 &0.950
		&0.945 &0.960 &0.953 &0.958 &0.950 &0.947\\
		&1  
		&1.000 &0.963
		&0.953 &0.960 &0.945 &0.947 &0.960 &0.955\\
		\hline
		\multirow{3}{*}{ACL}   &0 
		&$9.59\times10^6$ &2.643
		&5.258 &3.745 &3.063 &2.667 &2.406 &2.222\\
		&0.01  
		&$9.51\times10^6$ &2.640
		&5.222 &3.717 &3.057 &2.661 &2.398 &2.214\\
		&1  
		&$9.26\times10^6$ &2.614
		&5.267 &3.773 &3.068 &2.672 &2.414 &2.222\\
		\hline
		C.Time (s)  &&0.03 &29.02  &\multicolumn{6}{c}{2.166}\\
	\end{tabular}
\end{table}

\begin{table}
	\caption{\label{tab:p_600_ind} $N_b=624$, $b=12$, $p=600$, $s_0=10$, $\bSigma= \bI_p$. Performance on statistical inference. ``MLE" is the offline estimator obtained by fitting the traditional GLM, ``offline" corresponds to the offline debiased $\ell_1$-norm penalized estimator, and ``ODL" represents our proposed online debiased lasso estimator. Tuning parameter $\lambda$ is chosen from $T_\lambda= \{10^{-4}, 10^{-3}, 0.01, 0.05\}$. Simulation results are summarized over 200 replications. In the table, we report the $\lambda$ selected with highest frequency among 200 replications.}
	\centering
	\begin{tabular}{l l| r| r | rrrrrr}
		&$\beta_{0,r}$ &MLE  &offline &\multicolumn{6}{c}{ODL} \\
		data batch index &&&&2  &4 &6 &8  &10 &12\\
		$\lambda$ &&& &0.05  &0.01  &0.01  &0.01 &$10^{-4}$ &$10^{-4}$\\
		\hline
		\hline
		\multirow{3}{*}{A.bias} &0 
		&1.101 &0.014
		&0.034 &0.024 &0.020 &0.018 &0.015 &0.014\\
		&0.01 
		&1.553  &0.011
		&0.038 &0.014 &0.012 &0.009 &0.009 &0.010\\
		&1 
		&11.99  &0.069
		&0.168 &0.152 &0.152 &0.141 &0.131 &0.121\\
		\hline
		\multirow{3}{*}{ASE}  &0 
		&$1.82\times10^6$ &0.259
		&0.573 &0.412 &0.341 &0.297 &0.268 &0.245\\
		&0.01  
		&$1.80\times10^6$ &0.259
		&0.572 &0.411 &0.340 &0.297 &0.268 &0.246\\
		&1    
		&$1.81\times10^6$ &0.261
		&0.573 &0.411 &0.341 &0.297 &0.268 &0.246 \\
		\hline
		\multirow{3}{*}{ESE}  &0 
		&19.94 &0.254
		&0.578 &0.041 &0.034 &0.030 &0.028 &0.025\\
		&0.01 
		&19.72 &0.257
		&0.572 &0.397 &0.325 &0.287 &0.264 &0.251\\
		&1  
		&20.37 &0.262
		&0.579 &0.407 &0.334 &0.297 &0.262 &0.244\\
		\hline
		\multirow{3}{*}{CP}   &0 
		&1.000 &0.955
		&0.947 &0.947 &0.949 &0.948 &0.950 &0.950\\
		&0.01  
		&1.000 &0.953
		&0.945 &0.959 &0.954 &0.950 &0.956 &0.947\\
		&1  
		&1.000 &0.943
		&0.933 &0.935 &0.933 &0.925 &0.927 &0.928\\
		\hline
		\multirow{3}{*}{ACL}   &0 
		&$7.15\times10^6$ &1.016
		&2.246 &1.613 &1.335 &1.163 &1.050 &0.962\\
		&0.01  
		&$7.07\times10^6$ &1.016
		&2.241 &1.610 &1.334 &1.164 &1.050 &0.962\\
		&1  
		&$7.09\times10^6$ &1.022
		&2.246 &1.613 &1.335 &1.164 &1.050 &0.963\\
		\hline
		C.Time (min)  &&0.05 &117.5  &\multicolumn{6}{c}{7.08}\\
	\end{tabular}
\end{table}

\begin{table}
	\caption{\label{tab:p_600_ar1} $N_b=624$, $b=12$, $p=600$, $s_0=10$, $\bSigma= \{0.5^{|i-j|} \}_{i,j=1,\dots,p}$. Performance on statistical inference. ``MLE" is the offline estimator obtained by fitting the traditional GLM, ``offline" corresponds to the offline debiased $\ell_1$-norm penalized estimator, and ``ODL" represents our proposed online debiased lasso estimator. Tuning parameter $\lambda$ is chosen from $T_\lambda= \{10^{-4}, 10^{-3}, 0.01, 0.05\}$. Simulation results are summarized over 200 replications. In the table, we report the $\lambda$ selected with highest frequency among 200 replications.}
	\centering
	\begin{tabular}{l l| r| r | rrrrrr}
		&$\beta_{0,r}$ &MLE  &offline &\multicolumn{6}{c}{ODL} \\
		data batch index &&&&2  &4 &6 &8  &10 &12\\
		$\lambda$ &&& &0.05  &0.05  &$10^{-4}$  &$10^{-4}$ &$10^{-4}$  &$10^{-4}$\\
		\hline
		\hline
		\multirow{3}{*}{A.bias} &0 
		&1.594 &0.018
		&0.042 &0.034 &0.031 &0.029 &0.027 &0.025\\
		&0.01 
		&1.554 &0.016
		&0.044 &0.031 &0.026 &0.012 &0.010 &0.010\\
		&1 
		&11.85 &0.077
		&0.192 &0.168 &0.159 &0.147 &0.137 &0.134\\
		\hline
		\multirow{3}{*}{ASE}  &0 
		&$2.43\times10^6$ &0.307
		&0.576 &0.416 &0.347 &0.304 &0.274 &0.251\\
		&0.01  
		&$2.45\times10^6$ &0.307
		&0.577 &0.416 &0.347 &0.304 &0.274 &0.251\\
		&1    
		&$2.42\times10^6$ &0.308
		&0.576 &0.416 &0.347 &0.304 &0.275 &0.251\\
		\hline
		\multirow{3}{*}{ESE}  &0 
		&27.18 &0.300
		&0.576 &0.411 &0.341 &0.298 &0.269 &0.246\\
		&0.01 
		&27.68 &0.296
		&0.576 &0.403 &0.330 &0.287 &0.266 &0.240\\
		&1  
		&27.49 &0.291
		&0.559 &0.395 &0.326 &0.295 &0.257 &0.240\\
		\hline
		\multirow{3}{*}{CP}   &0 
		&1.000 &0.955
		&0.947 &0.949 &0.949 &0.949 &0.948 &0.947 \\
		&0.01  
		&1.000 &0.957
		&0.949 &0.958 &0.951 &0.958 &0.951 &0.956\\
		&1  
		&1.000 &0.951
		&0.946 &0.939 &0.929 &0.925 &0.931 &0.926\\
		\hline
		\multirow{3}{*}{ACL}   &0 
		&$9.51\times10^6$ &1.202
		&2.259 &1.631 &1.361 &1.193 &1.076 &0.982\\
		&0.01  
		&$9.59\times10^6$ &1.203
		&2.261 &1.632 &1.361 &1.193 &1.076 &0.983\\
		&1  
		&$9.48\times10^6$ &1.206
		&2.257 &1.629 &1.361 &1.193 &1.077 &0.983\\
		\hline
		C.Time (min)  &&0.06 &112.3  &\multicolumn{6}{c}{6.99}\\
	\end{tabular}
\end{table}

\spacingset{1.9}

\section{Real data analysis}\label{sec:data}
Graduated driver licensing programs are designed to restrict the driving privileges of new drivers. They typically include restrictions such as nighttime, expressway, and unsupervised driving. To monitor the effectiveness of such systems in reducing fatal teen crashes, we apply our proposed ODL procedure to analyze streaming data from the National Automotive Sampling System-Crashworthiness Data System (NASS CDS). Our primary interest was to evaluate the effectiveness of
graduated driver licensing, which is a nationwide legislature on novice drivers of age 21 or younger under various conditions of vehicle operation. In contrast, there are no restrictions on vehicle operation for older drivers (say, older than 65) in the current policy. To assess the effect of driver's age on driving safety, we compared age groups with respect to the risk of fatal crash when an accident happened. We first categorized the ``Age" variable into three groups: ``Age$<$21" representing the young group under a restricted graduate driver licensing, and ``Age$\geq$65" for the older group with a regular full driver's license, and those in between were treated as the reference group. Extent of ``Injury" in a crash is a binary outcome of interest, 1 for a moderate or severe injury, and 0 for minor or no injury. This outcome variable was created from the variable of Maximum Known Occupant Ais (MAIS), which indicates the single most severe injury level reported for each occupant. Other potential risk factors were also included in the model, including seat belt use, alcohol, speed limit, vehicle weight, air bag system deployed, drug involvement in this accident. Some interaction terms are also included such as driver's age group and alcohol, age group and sex etc.

Streaming data were formed by monthly accident data from the period of 7 years from January, 2009 to December, 2015, with $b=84$ data batches and a total of $N_b=2118$ crashed vehicles with complete records on both outcome and predictors. The number of predictors is $p=62$ including an intercept, 46 main effects and 15 interaction terms. Note that the first batch size $n_1=38<p=62$, and the smallest batch only contains 6 samples. We applied our proposed method to sequentially updated parameter estimates and standard errors for the regression coefficients in a penalized logistic regression model.

\begin{figure}[h]
	\centering
	\includegraphics[width=\linewidth]{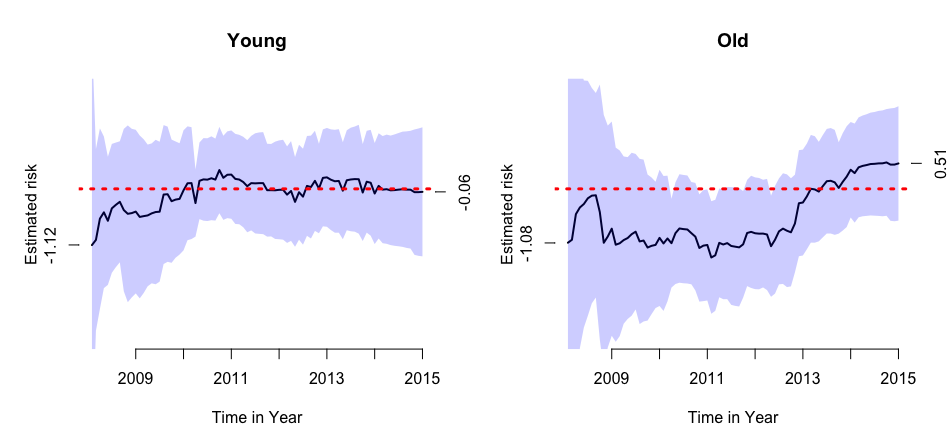}
	\caption{\label{fig:age_trace}Trace plot of the estimated risk and 95\% confidence bands of regression coefficients corresponding to young and old age groups, respectively. Numbers on two sides denote the estimated regression coefficients after the arrival of the first and last batches. The dashed line is the ``0" reference line.}
\end{figure}

As shown in Figure~\ref{fig:age_trace}, the 95\% pointwise confidence bands over the 84 batches became narrower for the old group but not the young one. Such phenomenon clearly differentiate the significance levels of these two age groups with respect to the severity of injury. Starting at almost the same estimated risk, the young age group gradually moves up and stay around 0 over the 84-month period. The trace plot ends at $-0.06$, meaning that the young age group has a slightly lower adjusted odds of fatal crash in comparison to the middle age group. This finding is consistent with the reported results in the literature that GDL is an effective policy to protect novice drivers from fatal injuries~\citep{GDL2014}. In contrast, the trace plot for the older age group shows an upward trend at around year 2013 and get stabilized at a positive value 0.51 since 2014. This observation indicates that the adjusted odds of fatality in a vehicle crash for the older group becomes significantly higher than the middle age group as more data accumulates over time. This may suggest a need on policy modification on restrictive vehicle operation for old drivers.

\begin{figure}[h]
	\centering
	\includegraphics[width=0.99\linewidth]{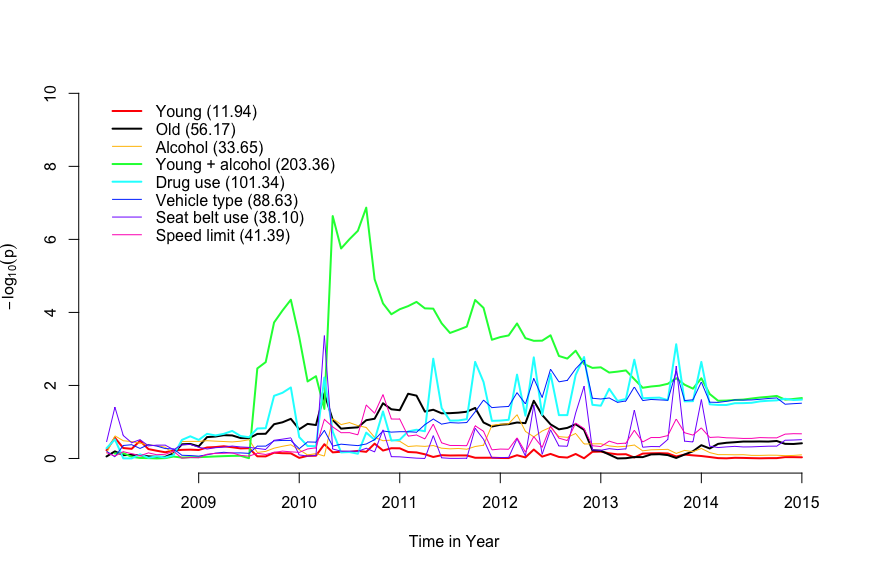}
	\caption{\label{fig:p_trace}Trace plot of $-\log_{10}p$ over monthly data from January, 2009 to December, 2015, each for one risk factor. Numbers of the $y$-axis are the negative logarithm $p$-values obtained by the Wald test and labels on the $x$-axis correspond to the last month of each year. The values in the brackets next to covariate names denote the areas under the $p$-value curves.}
\end{figure}

Figure~\ref{fig:p_trace} shows the trends of $-\log_{10}p$, $p$-values of the Wald test constructed using the online debiased lasso estimator and its standard error in the 10-base logarithm, for nine of the risk factors over 84 months. The trajectory of each covariate effect estimate show an evidence against the null $H_0:\beta^0_j=0$ as data accumulates over time. The interaction term between ``Age" and ``Alcohol" turns out to have the strongest association to the odds of fatality in a crash among all covariates shown in this plot. This is an overwhelming confirmation to the findings that driving under the influence of alcohol is significantly associated with the fatal crash even with the existence of graduate driver licensing to protect young drivers. This also explains why the main effect of young age group does not show a significant association with the odds of injury while the interaction term does. Moreover, some spikes appeared in the trajectories of most risk factors between year 2010 and 2014, but they get smoothed and stabilized thereafter. This might be related to small batch sizes or unbalanced categorical values in some monthly data batches. In addition, to characterize the overall significance level for each risk factor over the 84-month period, we propose to calculate a summary statistic as the area under the $p$-value curve. These values are included in the brackets after the name of each covariate in Figure~\ref{fig:p_trace}. They provide a clear summary of the overall significance level in existence of those spikes, and most of them are well aligned with the ranking of $p$-values at the end of year 2015.

Applying the proposed online debiased lasso procedure to the above CDS data analysis enabled us to visualize time-course patterns of data evidence accrual as well as stability and reproducibility of inference. As shown clearly in both Figures~\ref{fig:age_trace} and~\ref{fig:p_trace}, at the early stage of data streams, due to limited sample sizes and possibly sampling bias, both parameter estimates and test power may be unstable and even misleading. These potential shortcomings can be convincingly overcome when estimates and inferential quantities were continuously updated along with data streams, which eventually reached stability and reliable conclusions.

\section{Discussion}\label{sec:discussion}
In this paper, we propose an online debiased lasso method for statistical inference in high-dimensional GLMs. The method is applicable to 
streaming data, that is, only the historical summary statistics, not the raw historical data, are needed in updating the estimate at the current stage.
Under regularity conditions similar to those in the offline setting and mild conditions on the batch sizes, we prove the online debiased lasso (with an online correction term) is asymptotically normal.
The numerical studies further demonstrate the effectiveness of our algorithm
and support the theoretical results.

There are still many open questions in the area of online inference. First, the method we develop focuses on homogeneous data, where the streaming data is assumed to be i.i.d. sampled. It will be interesting to determine if we can obtain similar results under the non-homogeneous setting. Second, the loss function we consider in this paper is from the negative log-likelihood function. It is unclear whether other loss functions, including  nonsmooth robust loss functions such as the 
Huber's loss \citep{huber1964}, could be used in the online inference. Third, we did not address the issue about the online variable selection. The major difficulty in this problem is how to recover the significant variables which may be dropped at the early stages of the stream. We hope to address these interesting questions in the future.

\section*{Appendix}
The appendix contains the proofs of Theorems \ref{thm: lasso consistency}-\ref{thm: lasso asymptotic normality}.
\begin{proof}[Proof of Theorem \ref{thm: lasso consistency}]
	For the prior data batch $ \mathcal{D}_1 $, we have $ \widehat{\bbeta}^{(1)} = \bar{\bbeta}^{(1)} $ where $ \bar{\bbeta}^{(1)} $ is the offline lasso estimator defined in \eqref{lasso_offline}. Since the consistency of $ \bar{\bbeta}^{(1)} $ is well-established in \cite{vandegeer2014}, \eqref{lasso_consistency} holds when $ b = 1. $ Now we prove the
	consistency of $ \widehat{{\bbeta}}^{(b)} $ for an arbitrary  $ b \geq 2  $ by the mathematical induction.
	
	Suppose that $ {\widehat{\bbeta}}^{(b-1)} $ satisfies the statement in \eqref{lasso_consistency} with constant $ c_{1}^{(b-1)} $ and $ c_{2}^{(b-1)} $.  We claim that $ \lVert {\widehat{\bbeta}}^{(b)} - {{\bbeta}}^{0} \lVert_1 \leq c_{1}^{(b)} s_0 \lambda_b. $ Otherwise, we consider the following linear combination,
	\begin{eqnarray}\label{linear_combination}
		\widetilde{\bbeta}^{(b)} = t \widehat{\bbeta}^{(b)} + (1-t)\bbeta^{0}, \ \text{where} \ t = \frac{c_{1}^{(b)} s_0 \lambda_b}{c_{1}^{(b)} s_0 \lambda_b + \lVert {\widehat{\bbeta}}^{(b)} - {{\bbeta}}^{0} \lVert_1}.
	\end{eqnarray}
	Then $ \lVert {\widetilde{\bbeta}}^{(b)} - {{\bbeta}}^{0} \lVert_1 \leq c_{1}^{(b)} s_0 \lambda_b. $ Note that $ \lVert {\widetilde{\bbeta}}^{(b)} - {{\bbeta}}^{0} \lVert_1 \leq c_{1}^{(b)} s_0 \lambda_b/2 $ if and only if  $ \lVert {\widehat{\bbeta}}^{(b)} - {{\bbeta}}^{0} \lVert_1 \leq c_{1}^{(b)} s_0 \lambda_b. $ Therefore, it suffices to show $ \lVert {\widetilde{\bbeta}}^{(b)} - {{\bbeta}}^{0} \lVert_1 \leq c_{1}^{(b)} s_0 \lambda_b/2 $.
	
	Let $ \mL(\bbeta; \mD_j) = \sum_{i\in\mathcal{D}_j} d\left(y_i; g(\bx_i^\top\bbeta)\right) $, $ j = 1, \ldots, b $. Since the objective function defined in \eqref{lasso_step_j} is the convex function, we have
	\begin{eqnarray}
		&&\frac{1}{2N_b}\left\{ \mathcal{L}(\widetilde{\bbeta}^{(b)}; \mD_b) + \frac{1}{2}(\widetilde{\bbeta}^{(b)} - \widehat{\bbeta}^{(b-1)})^\top\widetilde{\bJ}^{(b-1)} (\widetilde{\bbeta}^{(b)} - \widehat{\bbeta}^{(b-1)})    \right\} + \lambda_{b}\lVert \widetilde{\bbeta}^{(b)} \lVert_1 \nonumber\\
		&\leq &\frac{1}{2N_b} \left\{ \mathcal{L}({\bbeta}^0; \mD_b) + \frac{1}{2}({\bbeta}^0 - \widehat{\bbeta}^{(b-1)})^\top \widetilde{\bJ}^{(b-1)} ({\bbeta}^0 - \widehat{\bbeta}^{(b-1)})    \right\} + \lambda_{b}\lVert \bbeta^0 \lVert_1, \label{lasso_inequality}
	\end{eqnarray}
	where $\widetilde{\bJ}^{(b-1)} = \sum_{j=1}^{b-1}\bJ^{(j)}(\widehat{\bbeta}^{(j)})$.
	Recall that $ \mathcal{D}_{b-1}^\star = \{\mathcal{D}_{1}, \ldots, \mathcal{D}_{b-1}\} $. A Taylor's expansion gives that
	\begin{align*}
		& \mL(\widetilde{\bbeta}^{(b)}; \mD_{b-1}^\star) - \mL({\bbeta}^0; \mD_{b-1}^\star) \\
		= & ~ \bU^{(b-1)}({\bbeta}^0) (\widetilde{\bbeta}^{(b)}- {\bbeta}^0) +  \frac{1}{2}(\widetilde{\bbeta}^{(b)}- {\bbeta}^0)^\top\left\{\widetilde{\bJ}^{(b-1)}(\bm{\xi})\right\} (\widetilde{\bbeta}^{(b)}- {\bbeta}^0),
	\end{align*}
	where $ \widetilde{\bJ}^{(b-1)}(\bm{\xi}) = \sum_{j=1}^{b-1}\bJ^{(j)}(\bm{\xi})$ and $ \bm{\xi} = t_2{\bbeta}^0 +(1-t_2)\widetilde{\bbeta}^{(b)} $ for some $ 0 < t_2 < 1. $
	Then
	\begin{eqnarray*}
		&&\frac{1}{2}(\widetilde{\bbeta}^{(b)} - \widehat{\bbeta}^{(b-1)})^\top\widetilde{\bJ}^{(b-1)} (\widetilde{\bbeta}^{(b)} - \widehat{\bbeta}^{(b-1)}) -  \frac{1}{2}({\bbeta}^0 - \widehat{\bbeta}^{(b-1)})^\top \widetilde{\bJ}^{(b-1)} ({\bbeta}^0 - \widehat{\bbeta}^{(b-1)}) \\
		&=&  \frac{1}{2}(\widetilde{\bbeta}^{(b)} - {\bbeta}^0)^\top \widetilde{\bJ}^{(b-1)} (\widetilde{\bbeta}^{(b)} - {\bbeta}^0) -  (\widetilde{\bbeta}^{(b)} - {\bbeta}^0)^\top \widetilde{\bJ}^{(b-1)} (\widehat{\bbeta}^{(b-1)} - {\bbeta}^0)
		\\
		& = & \mL(\widetilde{\bbeta}^{(b)}; \mD_{b-1}^\star) - \mL({\bbeta}^0; \mD_{b-1}^\star)  + \frac{1}{2}\sum_{j=1}^{b-1}(\widetilde{\bbeta}^{(b)} - {\bbeta}^0)^\top\left\{{\bJ}^{(j)}(\widehat{\bbeta}^{(j)}) - {\bJ}^{(j)}(\bm{\xi}) \right\} (\widetilde{\bbeta}^{(b)}- {\bbeta}^0) \\
		&& -\sum_{j=1}^{b-1}(\widetilde{\bbeta}^{(b)} - {\bbeta}^0)^\top \left\{{\bJ}^{(j)}(\widehat{\bbeta}^{(j)})\right\} (\widetilde{\bbeta}^{(b-1)} - {\bbeta}^0)
		-  \bU^{(b-1)}({\bbeta}^0) (\widetilde{\bbeta}^{(b)}- {\bbeta}^0) \\
		&:=& \mL(\widetilde{\bbeta}^{(b)}; \mD_{b-1}^\star) - \mL({\bbeta}^0; \mD_{b-1}^\star) + \Delta_1^{(b)} - \Delta_2^{(b)} - \Delta_3^{(b)}.
	\end{eqnarray*}
	Taking the above equation back into \eqref{lasso_inequality}, we obtain that
	\begin{eqnarray}\label{new lasso inequality}
		\frac{1}{2N_b} \mathcal{L}(\widetilde{\bbeta}^{(b)}; \mD_{b}^\star)      + \lambda_{b}\lVert \widetilde{\bbeta}^{(b)} \lVert_1 + \frac{1}{2N_b}(\Delta_1^{(b)} - \Delta_2^{(b)} - \Delta_3^{(b)})
		\leq \frac{1}{2N_b} \mathcal{L}({\bbeta}^0; \mD_b) + \lambda_{b}\lVert \bbeta^0 \lVert_1.
	\end{eqnarray}
	We introduce two notations,
	\begin{eqnarray*}
		\mathcal{E}({\bbeta})= \frac{1}{2N_b}\mathbb{E}\left\{ \mathcal{L}({\bbeta}; \mD_b^\star) - \mathcal{L}({\bbeta}^0; \mD_b^\star)  \right\}, ~
		v(\bbeta; \mD_b^\star) =  \mathcal{L}({\bbeta}; \mD_b^\star)  - \mathbb{E}\left\{ \mathcal{L}({\bbeta}; \mD_b^\star) \right\}, \ \bbeta \in \mathbb{R}^p,
	\end{eqnarray*}
	as the excess risk and the empirical process. Note that $ \mathcal{E}({\bbeta}) $ does not depend on $\mD_b^\star $ since the data is i.i.d. sampled. With the above notation, \eqref{new lasso inequality} could be further written as
	\begin{align*}
		\mathcal{E}(\widetilde{\bbeta}^{(b)}) + \lambda_{b}\lVert \widetilde{\bbeta}^{(b)} \lVert_1 &\leq -\frac{1}{2N_b}\left\{ v(\widetilde{\bbeta}^{(b)}; \mD_b^\star) - v({\bbeta}^0; \mD_b^\star)\right\} + \lambda_{b}\lVert \bbeta^0 \lVert_1
		+\frac{1}{2N_b}(\Delta_1^{(b)} - \Delta_2^{(b)} - \Delta_3^{(b)})\\
		&:= -\frac{1}{2N_b} \Delta_4^{(b)}  + \lambda_{b}\lVert \bbeta^0 \lVert_1
		+\frac{1}{2N_b}(\Delta_1^{(b)} - \Delta_2^{(b)} - \Delta_3^{(b)})
	\end{align*}
	Recall that $ \bX^{(b-1)}_\star = ((\bX^{(1)})^\top, \ldots, (\bX^{(b-1)})^\top)^\top  \in \mathbb{R}^{N_{b-1}\times p}$. The next lemma provides the upper bound of   $|\Delta_i^{(b)}|, i = 1, 2, 3, 4,$ whose proof is given at the end of Appendix.
	
	\vspace{0.2cm}
	
	\begin{lemma}\label{lemma 1 in consistency}
		Under the conditions of Theorem \ref{thm: lasso consistency}, with probability at least $ 1 - p^{-3}, $
		\begin{eqnarray*}
			&&|\Delta_1^{(b)}| \leq 2KLc^{(1)}_1s_0\lambda_1 \lVert \bX_\star^{(b-1)}(\widetilde{\bbeta}^{(b)}- {\bbeta}^0) \lVert^2_2,  \\
			&&| \Delta_2^{(b)}| \leq  K_2 \lVert  \bX_\star^{(b-1)}(\widetilde{\bbeta}^{(b)} - {\bbeta}^0) \lVert_2 \lVert  \bX_\star^{(b-1)}(\widetilde{\bbeta}^{(b-1)} - {\bbeta}^0) \lVert_2,\\
			&&| \Delta_3^{(b)}| \leq  {\lambda_{b-1}{N_{b-1}}} \lVert \widetilde{\bbeta}^{(b)} - {\bbeta}^0 \lVert_1/8, \ | \Delta_4^{(b)}| \leq {\lambda_{b}{N_{b}}}  \lVert \widetilde{\bbeta}^{(b)} - {\bbeta}^0 \lVert_1/8,
		\end{eqnarray*}
		where $ L $ is Lipschitz constant defined in Assumption \ref{assumption_1}, $ K = \lVert X \lVert_\infty $ and $ K_2 = \sup_{i \in \mathcal{D}_b^\star} \lvert {g'(\bx_i \bbeta)} \lvert$.
	\end{lemma}
	Note that $  {\lambda_{b-1}\sqrt{N_{b-1}}} = {\lambda_{b}\sqrt{N_{b}}} $ and the upper bound of $|\Delta_1|$ could be absorbed in the upper bound of $|\Delta_2|$. According to Lemma \ref{lemma 1 in consistency},
	\begin{eqnarray*}
		&&\Delta_1^{(b)} - \Delta_2^{(b)} - \Delta_3^{(b)} - \Delta_4^{(b)}\\
		&\leq & 2K_2 \lVert  \bX_\star^{(b-1)}(\widetilde{\bbeta}^{(b)} - {\bbeta}^0) \lVert_2 \lVert  \bX_\star^{(b-1)}(\widetilde{\bbeta}^{(b-1)} - {\bbeta}^0) \lVert_2 + \frac{1}{4}{\lambda_{b}{N_{b}}} \lVert \widetilde{\bbeta}^{(b)} - {\bbeta}^0 \lVert_1\\
		&\leq& 2K_2\sqrt{N_{b-1}c_2^{(b-1)}s_0}\lambda_{b-1}\lVert  \bX_\star^{(b)}(\widetilde{\bbeta}^{(b)} - {\bbeta}^0) \lVert_2  + \frac{1}{4}c_1^{(b)}{{N_{b}}}s_0\lambda_{b}^2.
	\end{eqnarray*}
	Consequently, 
	\begin{eqnarray}
		&&\mathcal{E}(\widetilde{\bbeta}^{(b)})  + \lambda_{b}\lVert \widetilde{\bbeta}^{(b)} \lVert_1  - \lambda_{b}\lVert \bbeta^0 \lVert_1 \nonumber\\
		&\leq& \frac{1}{2N_b} \left\{ 2K_2\sqrt{N_{b-1}c_2^{(b-1)}s_0}\lambda_{b-1}\lVert  \bX_\star^{(b)}(\widetilde{\bbeta}^{(b)} - {\bbeta}^0) \lVert_2 +  \frac{1}{4}c_1^{(b)}{{N_{b}}}s_0\lambda_{b}^2 \right\} \nonumber\\
		&=& \left\{ K_2\sqrt{\frac{c_2^{(b-1)}s_0}{N_b}}\lVert  \bX_\star^{(b)}(\widetilde{\bbeta}^{(b)} - {\bbeta}^0) \lVert_2 +  \frac{1}{8}c_1^{(b)}s_0\lambda_{b}\right\}\lambda_{b}  := \Delta^{(b)}\lambda_{b}.\label{lasso: basic inequality 2}
	\end{eqnarray}
	The left part follows the standard argument. Recall that $ S_0 = \{j: \beta_{j}^0 \neq 0\} $. For $ \bbeta \in \mathbb{R}^p$, we define that $ \bbeta_{S_0} := (\beta_{j,S_0})_{j=1}^p $ where $ \beta_{j,S_0} = \beta_j \mathbf{1}\{j \in S_0\}. $ Then, $ \bbeta =  \bbeta_{S_0} + \bbeta_{S^c_0}.$ From \eqref{lasso: basic inequality 2}, we have
	\begin{eqnarray*}
		\mathcal{E}(\widetilde{\bbeta}^{(b)})   + \lambda_{b}\lVert \widetilde{\bbeta}^{(b)}_{S^c_0}\lVert_1  \leq  \lambda_{b}\lVert \bbeta^0 \lVert_1 -  \lambda_{b}\lVert \widetilde{\bbeta}^{(b)}_{S_0} \lVert_1 + \Delta^{(b)}\lambda_{b}
		\leq  \lambda_{b}\lVert \bbeta^0  -   \widetilde{\bbeta}^{(b)}_{S_0} \lVert_1 + \Delta^{(b)}\lambda_{b}.
	\end{eqnarray*}
	Here we need a discussion about the value of $ \lVert \bbeta^0  -   \widetilde{\bbeta}^{(b)}_{S_0} \lVert_1 $.\\
	\textbf{Case I.} We suppose that
	$ \lVert \bbeta^0  -   \widetilde{\bbeta}^{(b)}_{S_0} \lVert_1 \geq \Delta^{(b)}/2
	$. Then,
	\begin{equation*}
		\mathcal{E}(\widetilde{\bbeta}^{(b)})   + {\lambda_{b}}\lVert \widetilde{\bbeta}^{(b)}_{S^c_0}\lVert_1
		\leq  {3\lambda_{b}}\lVert \bbeta^0  -   \widetilde{\bbeta}^{(b)}_{S_0} \lVert_1,
	\end{equation*}
	which implies $ \lVert \widetilde{\bbeta}^{(b)}_{S^c_0}\lVert_1 \leq 3\lVert \bbeta^0  -   \widetilde{\bbeta}^{(b)}_{S_0} \lVert_1 $. Then, we can adopt the {empirical compatibility condition, that is,
		\begin{equation*}
			\lVert \bbeta^0  -   \widetilde{\bbeta}^{(b)}_{S_0} \lVert_1^2 \leq \frac{s_0}{\phi_0^2N_b} \lVert \bX^{(b)}_\star(\bbeta^0  -   \widetilde{\bbeta}^{(b)}) \lVert^2_2,
		\end{equation*}
		where $ \phi_0 > 0 $ is the compatibility constant.} Thus,
	\begin{eqnarray*}
		&&\mathcal{E}(\widetilde{\bbeta}^{(b)})   + {\lambda_{b}}\lVert \widetilde{\bbeta}^{(b)}_{S^c_0}\lVert_1 + {\lambda_{b}}\lVert \bbeta^0  -   \widetilde{\bbeta}^{(b)}_{S_0} \lVert_1
		\leq  \frac{2\lambda_{b} }{\phi_0}\sqrt{\frac{s_0}{N_b}}\lVert \bX^{(b)}_\star(\bbeta^0  -   \widetilde{\bbeta}^{(b)}) \lVert_2  + \Delta^{(b)}\lambda_{b}\\
		&=& \left(\frac{2}{\phi_0} + K_2\sqrt{c_2^{(b-1)}} \right)
		\lambda_{b} \sqrt{\frac{s_0}{N_b}}\lVert \bX^{(b)}_\star(\bbeta^0  -   \widetilde{\bbeta}^{(b)}) \lVert_2 + \frac{1}{8}c_1^{(b)}s_0\lambda_{b}^2\\
		&:=& C\lambda_{b}\sqrt{\frac{s_0}{N_b}}\lVert \bX^{(b)}_\star(\bbeta^0  -   \widetilde{\bbeta}^{(b)}) \lVert_2 + \frac{1}{8}c_1^{(b)}s_0\lambda_{b}^2
	\end{eqnarray*}
	Let $ k_2 = 1/ \sup_{i \in \mathcal{D}_b^\star} \lvert {4/g'(\bx_i \bbeta)} \lvert. $ Based on Lemma \ref{lemma 1 in consistency}, we have
	\begin{equation*}
		\frac{k_2}{N_b} \lVert \bX^{(b)}_\star(\bbeta^0  -   \widetilde{\bbeta}^{(b)}) \lVert^2_2 \leq \mathcal{E}(\widetilde{\bbeta}^{(b)}) +  \frac{1}{2N_b} \Delta_4^{(b)} \leq \mathcal{E}(\widetilde{\bbeta}^{(b)}) +  \frac{1}{16}c_1^{(b)}s_0\lambda_{b}^2,
	\end{equation*}
	which is also known as the margin condition \citep{van2008high}.
	{ Next we apply the arithmetic mean-geometric mean inequality  and obtain}
	\begin{equation}\label{equation_1: case_1}
		C\lambda_{b}\sqrt{\frac{s_0}{N_b}}\lVert \bX^{(b)}_\star(\bbeta^0  -   \widetilde{\bbeta}^{(b)}) \lVert_2 \leq
		\frac{C^2}{2k_2}s_0\lambda_{b}^2  + \frac{\mathcal{E}(\widetilde{\bbeta}^{(b)})}{2} +  \frac{1}{32}c_1^{(b)}s_0\lambda_{b}^2.
	\end{equation}
	Then, it follows that 
	\begin{equation}\label{equation_2: case_1}
		\frac{\mathcal{E}(\widetilde{\bbeta}^{(b)})}{2}   + {\lambda_{b}}\lVert \bbeta^0  -   \widetilde{\bbeta}^{(b)} \lVert_1
		\leq    \left( \frac{C^2}{2k_2}+ \frac{5}{32}c_1^{(b)}\right) s_0\lambda_{b}^2 .
	\end{equation}
	On the one hand, since $  \mathcal{E}(\widetilde{\bbeta}^{(b)}) > 0$, we obtain that
	$$ \lVert \bbeta^0  -   \widetilde{\bbeta}^{(b)} \lVert_1 \leq \left( \frac{C^2}{2k_2}+ \frac{5}{32}c_1^{(b)}\right) s_0\lambda_{b}.$$
	With the suitable choice of $ c_2^{(b-1)} $,
	\begin{equation*}
		\frac{C^2}{2k_2}+ \frac{5}{32}c_1^{(b)} = \frac{1}{2k_2}\left(\frac{2}{\phi_0} + K_2\sqrt{c_2^{(b-1)}} \right)^2 + \frac{5}{32}c_1^{(b)} \leq \frac{13}{32} c_1^{(b)},
	\end{equation*}
	we have
	\begin{equation*}
		\lVert \bbeta^0  -   \widetilde{\bbeta}^{(b)} \lVert_1
		\leq    \frac{1}{2}c_1^{(b)} s_0\lambda_{b} .
	\end{equation*}
	{Here we require that $ K_2\sqrt{c_2^{(1)}} \geq 2/\phi_0$,  and $ 8K_2^2c_2^{(b-1)}/k_2 = c_1^{(b)}.$ }
	On the other hand, we combine \eqref{equation_1: case_1} and \eqref{equation_2: case_1}, and obtain
	\begin{equation*}
		\lVert \bX^{(b)}_\star(\bbeta^0  -   \widetilde{\bbeta}^{(b)}) \lVert_2  \leq   \frac{5C}{2k_2}\sqrt{s_0N_b}\lambda_{b}.
	\end{equation*}
	Again, since $ {C} \leq {2K_2\sqrt{c_2^{(b-1)}}} $, we obtain
	\begin{equation}\label{case_1_estimation_error}
		\lVert \bX^{(b)}_\star(\bbeta^0  -   \widetilde{\bbeta}^{(b)}) \lVert^2_2 \leq \left( \frac{5K_2}{k_2}\right)^2c_2^{(b-1)}s_0\lambda_{b}^2N_b \leq c_2^{(b)}s_0\lambda_{b}^2N_b.
	\end{equation}
	\textbf{Case II.} We suppose that
	$ \lVert \bbeta^0  -   \widetilde{\bbeta}^{(b)}_{S_0} \lVert_1 < \Delta^{(b)}/2$. Then,
	\begin{eqnarray*}
		&&\mathcal{E}(\widetilde{\bbeta}^{(b)})   + {\lambda_{b}}\lVert \widetilde{\bbeta}^{(b)}_{S^c_0}\lVert_1 + {\lambda_{b}}\lVert \bbeta^0  -   \widetilde{\bbeta}^{(b)}_{S_0} \lVert_1
		\leq  2 \Delta^{(b)}\lambda_{b} \\
		&=& 2K_2
		\lambda_{b} \sqrt{\frac{c_2^{(b-1)}s_0}{N_b}}\lVert \bX^{(b)}_\star(\bbeta^0  -   \widetilde{\bbeta}^{(b)}) \lVert_2 + \frac{1}{4}c_1^{(b)}s_0\lambda_{b}^2
	\end{eqnarray*}
	Then, by the margin condition, it is straightforward to show that
	\begin{equation}\label{case_2_estimation_error}
		\lVert \bbeta^0  -   \widetilde{\bbeta}^{(b)} \lVert_1 \leq \frac{1}{2}c^{(b)}_1s_0\lambda_{b}, \ \lVert \bX(\bbeta^0  -   \widetilde{\bbeta}^{(b)}) \lVert_2^2  \leq   \left( \frac{6K_2}{k_2}\right)^2c_2^{(b-1)}s_0\lambda_{b}^2N_b = c_2^{(b)}s_0\lambda_{b}^2N_b.
	\end{equation}
	Note that in both cases, we obtain $  \lVert \bbeta^0  -   \widetilde{\bbeta}^{(b)} \lVert_1 \leq c^{(b)}_1s_0\lambda_{b}/2. $ According to \eqref{linear_combination}, we show our claim $  \lVert \bbeta^0  -   \widehat{\bbeta}^{(b)} \lVert_1 \leq c^{(b)}_1s_0\lambda_{b}. $ The remaining step is to repeat the above argument and obtain the upper bound of the estimation error at \eqref{case_1_estimation_error} or \eqref{case_2_estimation_error}. It is worth pointing out that $  c_1^{(b)} = 8K_2^2c_2^{(b-1)}/k_2 $ and $ c_2^{(b)} = 36K_2^2c_2^{(b-1)}/k_2^2, $ where $ K_2 $ and $ k_2 $ does not depend on $ b. $ The proof is finished by taking a union bound on the events considered in Lemma \ref{lemma 1 in consistency}.
\end{proof}

Before proving Theorem \ref{thm: lasso consistency 2}, we cite Lemma 3 in \cite{han2021online} to compute the cumulative error.
\begin{lemma}\citep{han2021online}\label{lemma_2}
	Let $ n_j $ and $ N_j $ be the batch size and the cumulative batch size respectively when the $ j$-th data arrives,  $ j = 1,\ldots, b. $ Then,
	\begin{eqnarray}
		\sum_{j=1}^{b}\frac{n_j}{N_j} &\leq& 1 + \log \frac{N_b}{n_1}, \label{cumulative_1}\\
		\sum_{j=1}^b \frac{n_j}{\sqrt{N_j}} &\leq& 2\sqrt{N_b}. \label{cumulative_2}
	\end{eqnarray}
\end{lemma}

\begin{proof}[Proof of Theorem \ref{thm: lasso consistency 2}]
	Recall that
	\begin{equation*}
		\bm{\widehat{\gamma}}_r^{(j)}= \underset{\bm{\gamma}\in\mathbb{R}^{(p-1)}}{\arg\min}
		\left\{\frac{1}{2N_j}\left(\widetilde{\bJ}^{(j)}_{r,r} - 2 \widetilde{\bJ}^{(j)}_{r,-r}\bgamma +\bgamma^\top\widetilde{\bJ}_{-r,-r}^{(j)}\bgamma \right) +\lambda_j\|\bm{\gamma} \|_1
		\right\}
	\end{equation*}
	is a standard lasso estimator. The proof follows the standard argument as long as estimated information matrix $ \widetilde{\bJ}^{(j)}/N_j $ satisfies the compatibility condition. Therefore, it suffices to demonstrate the compatibility condition holds. Recall that $ \widetilde{\bJ}^{(j)}/N_j = \sum_{i=1}^{j} \bJ^{(i)}(\widehat{\bbeta}^{(i)})/N_j.$ Then we have
	\begin{eqnarray*}
		&&	\widetilde{\bJ}^{(j)}/N_j - \bJ^{0} = \frac{1}{N_j}\sum_{i=1}^{j}\left\{ \bJ^{(i)}(\widehat{\bbeta}^{(i)}) -  \mathbb{E}[\bJ^{(i)}({\bbeta}^{0})]\right\} \\
		&=& \frac{1}{N_j}\sum_{i=1}^{j}\left\{ \bJ^{(i)}(\widehat{\bbeta}^{(i)}) - \bJ^{(i)}({\bbeta}^{0})\right\} +\frac{1}{N_j}\sum_{i=1}^{j}\left[ \bJ^{(i)}({\bbeta}^{0}) - \mathbb{E}\{\bJ^{(i)}({\bbeta}^{0})\}\right].
	\end{eqnarray*}
	According to the error bound of $ \widehat{\bbeta}^{(i)} $ provided in Theorem \ref{thm: lasso consistency} and \eqref{cumulative_2}, we obtain
	\begin{equation*}
		\left\lVert \frac{1}{N_j}\sum_{i=1}^{j}\left\{ \bJ^{(i)}(\widehat{\bbeta}^{(i)}) - \bJ^{(i)}({\bbeta}^{0})\right\} \right\lVert_\infty \leq \frac{K^2}{N_j}\sum_{i=1}^{j} {n_i} \lVert {\widehat{\bbeta}}^{(i)} - {{\bbeta}}^{0} \lVert_1 = \mathcal{O}_p\left( c_1^{(j)}s_0\sqrt{\frac{\log(p)}{N_j}}\right),
	\end{equation*}
	where $ K $ is defined in \textit{(A1)} of Assumption \ref{assumption_1}.
	Meanwhile, based on the Hoeffding's inequality, it follows
	\begin{equation*}
		\left\lVert \frac{1}{N_j}\sum_{i=1}^{j}\left[ \bJ^{(i)}({\bbeta}^{0}) - \mathbb{E}\{\bJ^{(i)}({\bbeta}^{0})\}\right] \right\lVert_\infty  = \mathcal{O}_p\left(\sqrt{\frac{\log(p)}{N_j}}\right).
	\end{equation*}
	In summary,
	\begin{equation*}
		\left\lVert	\widetilde{\bJ}^{(j)}/N_j - \bJ^{0}\right\lVert_\infty =  \mathcal{O}_p\left( c_1^{(j)}s_0\sqrt{\frac{\log(p)}{N_j}}\right).
	\end{equation*}
	Consequently, the compatibility condition holds through Corollary 6.8 in \cite{buhlmann2011statistics}. 
\end{proof}

The remaining part is to establish the asymptotic normality of the online debiased lasso estimator.
\begin{proof}[Proof of Theorem \ref{thm: lasso asymptotic normality}]
	We first deal with the debiased term. Recall that $ \widetilde{\bgamma}^{(j)} \in \mathbb{R}^p $ is the extension of $ \widehat
	{\bgamma}^{(j)} $.
	Then
	\begin{eqnarray*}
		&&\sum_{j=1}^{b}\{\bm{\widehat{z}}^{(j)}_{r}\}^\top\left\{\by^{(j)} - g(\bX^{(j)}{\widehat{\bbeta}}^{(j)} )\right\} \\
		&=& \sum_{j=1}^{b} \{{\widehat{\bz}}^{(j)}_{r}\}^\top\left\{\by^{(j)} - g(\bX^{(j)}{{\bbeta}^0} )\right\} - \{\widetilde{\bgamma}^{(j)}_r\}^\top  \{\bX^{(j)}\}^\top\left\{ g(\bX^{(j)}{{\bbeta}^0}) - g(\bX^{(j)}{\widehat{\bbeta}}^{(j)} )\right\}
	\end{eqnarray*}	
	Note that for $ 1 \leq k \leq n, $ we let $ \bX^{(j)}_k \in \mathbb{R}^{1\times p} $ denote $ k$-th row in $ \bX^{(j)} $, namely the covariates of $ k$-th data in $ j$-th batch. According to  the mean value theorem,
	\begin{eqnarray*}
		g(\bX^{(j)}_k{\widehat{\bbeta}}^{(j)})  = g(\bX^{(j)}_k{{\bbeta}}^0) - g'( \eta_k^{(j)})\bX^{(j)}_k({\bbeta}^0 - {\widehat{\bbeta}}^{(j)}),
	\end{eqnarray*}
	where $ \eta_k^{(j)} \in [\bX^{(j)}_k{{\bbeta}^0}, \bX^{(j)}_k{\widehat{\bbeta}}^{(j)}].$
	Let $ \bm{\eta}^{(j)} = (\eta_1^{(j)}, \ldots,  \eta_{n_j}^{(j)})^\top $ and $\bm{\Lambda}^{(j)} \in \mathbb{R}^{n_j\times n_j}$ is diagonal matrix with the diagonal element $ \{g'(\bX^{(j)}_k{\widehat{\bbeta}}^{(j)}) - g'({\eta}_k^{(j)}) \}_{k=1}^{n_j} $. As a result,
	\begin{eqnarray*}
		&&\sum_{j=1}^{b}\{\bm{\widehat{z}}^{(j)}_{r}\}^\top\left\{\by^{(j)} - g(\bX^{(j)}{\widehat{\bbeta}}^{(j)} )\right\} \\
		&=& \sum_{j=1}^{b} \{{\widehat{\bz}}^{(j)}_{r}\}^\top\left\{\by^{(j)} - g(\bX^{(j)}{{\bbeta}}^0 )\right\} - \sum_{j=1}^{b}\{\widetilde{\bgamma}^{(j)}_r\}^\top \bJ^{(j)}(\widehat{\bbeta}^{(j)})({\bbeta}^0 - {\widehat{\bbeta}}^{(j)}) +\Pi_1,
	\end{eqnarray*}	 
	where $ \Pi_1 = \sum_{j=1}^{b} \{\bm{\widehat{z}}^{(j)}_{r}\}^\top \bm{\Lambda}^{(j)}\bX^{(j)}({\bbeta}^0 - {\widehat{\bbeta}}^{(j)}) $.
	Now, we focus on the online debiased lasso estimator defined in \eqref{online_debiased_algorithm}. According to the above results,
	\begin{eqnarray*}
		&&\widehat{\tau}_r^{(b)} ({\widehat{\beta}}^{(b)}_{\text{on}, r} - \beta_{r}^{0})\\
		&=&  \sum_{j=1}^{b}\{\bm{\widehat{z}}^{(j)}_{r}\}^\top\left\{\by^{(j)} - g(\bm{X}^{(j)}\bm{\widehat{\beta}}^{(j)} )\right\} + \sum_{j=1}^{b}\{\widetilde{\bgamma}^{(j)}_r\}^\top\bJ^{(j)}(\widehat{\bbeta}^{(j)})({\widehat{\bbeta}}^{(b)} - {\widehat{\bbeta}}^{(j)} )  + \widehat{\tau}_r^{(b)} (\widehat{{\beta}}^{(b)}_{r} - \beta_{r}^{0})\\
		&=& \sum_{j=1}^{b} \{{\widehat{\bz}}^{(j)}_{r}\}^\top\left\{\by^{(j)} - g(\bX^{(j)}{{\bbeta}^0} )\right\} +\Pi_1  + \sum_{j=1}^{b} \{\widetilde{\bgamma}^{(j)}_r\}^\top \bJ^{(j)}(\widehat{\bbeta}^{(j)})({\widehat{\bbeta}}^{(b)} - {\bbeta}^0)  + \widehat{\tau}_r^{(b)} (\widehat{{\beta}}^{(b)}_{r} - \beta_{r}^{0}) \\
		&=& \sum_{j=1}^{b} \{{\widehat{\bz}}^{(j)}_{r}\}^\top\left\{\by^{(j)} - g(\bX^{(j)}{{\bbeta}^0} )\right\} + \Pi_1+ \Pi_2 + \Pi_3,
	\end{eqnarray*}	
	where
	\begin{eqnarray*}
		\Pi_2  &=& \sum_{j=1}^{b} \{\widetilde{\bgamma}^{(b)}_r\}^\top\bJ^{(j)}(\widehat{\bbeta}^{(j)})  ({\widehat{\bbeta}}^{(b)} - {\bbeta}^0) -\sum_{j=1}^{b} \{\widetilde{\bgamma}^{(b)}_r\}^\top\bJ^{(j)}_r(\widehat{\bbeta}^{(j)})  (\widehat{{\beta}}^{(b)}_{r} - \beta_{r}^{0})  ,\\
		\Pi_3  &=& \sum_{j=1}^{b} \{\widetilde{\bgamma}^{(j)}_r\}^\top\bJ^{(j)}(\widehat{\bbeta}^{(j)})  ({\widehat{\bbeta}}^{(b)} - {\bbeta}^0)- \sum_{j=1}^{b} \{\widetilde{\bgamma}^{(b)}_r\}^\top\bJ^{(j)}(\widehat{\bbeta}^{(j)})  ({\widehat{\bbeta}}^{(b)} - {\bbeta}^0)  .
	\end{eqnarray*}	
	Therefore, the remaining part is to demonstrate that  $ |\Pi_i| = o_p(\sqrt{N_b}), i = 1, 2, 3. $
	
	First of all, according to the Lipschitz condition of $ g'(\cdot), $
	\begin{eqnarray*}
		\left\lvert g'({\eta}^{(j)}_k)\bX^{(j)}_k({\bbeta}^0 - {\widehat{\bbeta}}^{(j)}) - g'(\bX^{(j)}_k{\widehat{\bbeta}}^{(j)})\bX^{(j)}_k({\bbeta}^0 - {\widehat{\bbeta}}^{(j)}) \right\lvert
		\leq \left\{ \bX^{(j)}_k({\bbeta}^0 - {\widehat{\bbeta}}^{(j)})\right\}^2.
	\end{eqnarray*}
	Then,
	\begin{eqnarray*}
		\left\lVert\bm{\Lambda}^{(j)}\bX^{(j)}({\bbeta}^0 - {\widehat{\bbeta}}^{(j)})\right\lVert_1  \leq \left\lVert \bX^{(j)}({\bbeta}^0 - {\widehat{\bbeta}}^{(j)})\right\lVert_2^2 \leq K^2 (c_1^{(j)}s_0\lambda_{j})^2n_j,
	\end{eqnarray*}
	where the last inequality is from Theorem \ref{thm: lasso consistency} and the bounded assumption of $ \bX. $ Meanwhile, the boundedness of $ \lVert  \widehat{\bz}_r^{(j)}\lVert_\infty  $ could be shown by $(A2)$ in Assumption \ref{assumption_1} and Theorem \ref{thm: lasso consistency 2}.
	Therefore, we find an upper bound for $|\Pi_1|$:
	\begin{eqnarray*}
		|\Pi_1| &\leq& \sum_{j=1}^{b}\lVert  \widehat{\bz}_r^{(j)}\lVert_\infty \left\lVert\bm{\Lambda}^{(j)}\bX^{(j)}({\bbeta}^0 - {\widehat{\bbeta}}^{(j)})\right\lVert_1\\
		&\leq&\sum_{j=1}^{b} K^2 \lVert  \widehat{\bz}_r^{(j)}\lVert_\infty (c_1^{(j)}s_0\lambda_{j})^2n_j = \mathcal{O}_p \left(c_1^{(b)} s_0^2 \log (p)\log(N_b) \right),
	\end{eqnarray*}
	where the last equation is from \eqref{cumulative_1} in Lemma \ref{lemma_2}.
	
	Next, we adopt the Karush-Kuhn-Tucker (KKT) conditions. Let $ \be_r \in \mathbb{R}^{p} $ denote the zero-vector except that the $ r$-th element is one.
	\begin{eqnarray*}
		|\Pi_2|
		&=&  \left\lvert \left( \{\widetilde{\bgamma}^{(b)}_r\}^\top\widetilde{\bJ}^{(b)} - \{\widetilde{\bgamma}^{(b)}_r\}^\top\widetilde{\bJ}^{(b)}_r\be_r^\top  \right) ({\widehat{\bbeta}}^{(b)} - {\bbeta}^0 )  \right\lvert            \\
		&\leq& \left\lVert\{\widetilde{\bgamma}^{(b)}_r\}^\top\widetilde{\bJ}^{(b)} - \{\widetilde{\bgamma}^{(b)}_r\}^\top\widetilde{\bJ}^{(b)}_r\be_r^\top  \right\lVert_\infty \left\lVert{\widehat{\bbeta}}^{(b)} - {\bbeta}^0\right\lVert_1\\
		&\leq& N_b\lambda_b \times c^{(b)}_1 s_0\lambda_{b} = c^{(b)}_1 s_0N_b\lambda_{b}^2 = \mathcal{O}_p\left( c^{(b)}_1s_0\log(p)\right).
	\end{eqnarray*}
	Finally, we apply Theorem \ref{thm: lasso consistency 2} and \eqref{cumulative_2} in Lemma \ref{lemma_2}.
	\begin{eqnarray*}
		|\Pi_3| &=& \left\lvert \left\{\sum_{j=1}^{b} \{\widetilde{\bgamma}^{(b)}_r - \widetilde{\bgamma}^{(j)}_r\}^\top\bJ^{(j)}(\widehat{\bbeta}^{(j)})\right\}({\widehat{\bbeta}}^{(b)} - {\bbeta}^0) \right\lvert\\
		&\leq& \left\lVert{\widehat{\bbeta}}^{(b)} - {\bbeta}^0\right\lVert_1 \sum_{j=1}^{b} \left\lVert \widetilde{\bgamma}^{(b)}_r - \widetilde{\bgamma}^{(j)}_r\right\lVert_1 \left\lVert\bJ^{(j)}(\widehat{\bbeta}^{(j)})\right\lVert_\infty \\
		&\leq& c_1^{(b)}s_0\lambda_{b} \left( \sum_{j=1}^{b} c s_r\lambda_{j}n_j\right)  = \mathcal{O}_p\left(c_1^{(b)}s_0s_r\log(p) \right).
	\end{eqnarray*}	
	Based on the conditions in Theorem \ref{thm: lasso asymptotic normality}, we conclude  $ |\Pi_i| = o_p(\sqrt{N_b}), i = 1, 2, 3. $ In summary,
	\begin{eqnarray*}
		\widehat{\tau}_r^{(b)}({\widehat{\beta}}^{(b)}_{\text{on}, r} - \beta_{r}^{0}) =  \sum_{j=1}^{b} \{{\widehat{\bz}}^{(j)}_{r}\}^\top\left\{\by^{(j)} - g(\bX^{(j)}{{\bbeta}^0} )\right\} + o_p(\sqrt{N_b}),
	\end{eqnarray*}
	where we could further refer the terms in the right parts as $ W_r $ and $ V_r $ in Theorem \ref{thm: lasso asymptotic normality}.
	
\end{proof}
For the sake of completeness, we provide the proof of Lemma \ref{lemma 1 in consistency} at the end of Appendix.
\begin{proof}[Proof of Lemma \ref{lemma 1 in consistency}]
	
	We start from $ |\Delta_1^{(b)}|. $ Recall that
	\begin{align*}
		\Delta_1^{(b)} = \frac{1}{2}\sum_{j=1}^{b-1}(\widetilde{\bbeta}^{(b)} - {\bbeta}^0)^\top\left\{{\bJ}^{(j)}(\widehat{\bbeta}^{(j)}) - {\bJ}^{(j)}(\bm{\xi}) \right\} (\widetilde{\bbeta}^{(b)}- {\bbeta}^0)
	\end{align*}
	Let $ \bv^{(j)} = \bX^{(j)}(\widetilde{\bbeta}^{(b)}- {\bbeta}^0)$, $ \bw^{(j)} = g'(\bX^{(j)} \widehat{{\bbeta}}^{(j)}) -  g'(\bX^{(j)}\bm{\xi}) $ and $ \text{diag}(\bw^{(j)}) $ denote the diagonal matrix with diagonal element $ \bw^{(j)}. $ Then,
	\begin{align*}
		\Delta_1^{(b)} = \frac{1}{2} \sum_{j=1}^{b-1} (\bv^{(j)})^\top \text{diag}(\bw^{(j)})\bv^{(j)} = \frac{1}{2} \bv^\top \text{diag}(\bw)\bv,
	\end{align*}
	where $ \bv = ((\bv^{(1)})^\top, \ldots, (\bv^{(b-1)})^\top)^\top $ and $ \bw = ((\bw^{(1)})^\top, \ldots, (\bw^{(b-1)})^\top)^\top $. As a result,
	\begin{align*}
		|\Delta_1^{(b)}| \leq \frac{1}{2} \lVert \bv \lVert_2^2 \lVert \text{diag}(\bw) \lVert_2 = \frac{1}{2} \lVert \bX_\star^{(b-1)}(\widetilde{\bbeta}^{(b)}- {\bbeta}^0) \lVert^2_2 \lVert \bw \lVert_\infty.
	\end{align*}
	The left part is to find the upper bound of $ \lVert \bw \lVert_\infty $, that is,
	\begin{eqnarray*}
		\lVert \bw \lVert_\infty &=& \max_{1\leq j\leq b-1} \lVert \bw^{(j)} \lVert_\infty = \max_{1\leq j\leq b-1} \lVert g'(\bX^{(j)} \widehat{{\bbeta}}^{(j)}) -  g'(\bX^{(j)}\bm{\xi}) \lVert_\infty \\
		&\leq& L\max_{1\leq j\leq b-1} \lVert \bX^{(j)} (\widehat{{\bbeta}}^{(j)} -  \bm{\xi}) \lVert_\infty \leq 3KL\max_{1\leq j\leq b-1} c^{(j)}_1s_0\lambda_{j}
	\end{eqnarray*}
	For ease of understanding, we could assume that $ c^{(1)}_1s_0\lambda_{1} = \max_{1\leq j\leq b-1} c^{(j)}_1s_0\lambda_{j} $ since the first data batch containing the least information may lead to the largest estimation error. Besides that, this assumption will not affect the outcome. The upper bound of $|\Delta_1|$ could always be absorbed in the upper bound of $|\Delta_2|$ due to $ b = o(\log N_b) $.
	
	For $ \Delta_2 $, the proof is similar as the above one and we omit it.
	
	Recall that $ |\Delta_3 ^{(b)}| =  \bU^{(b-1)}({\bbeta}^0) (\widetilde{\bbeta}^{(b)}- {\bbeta}^0)$. Note that $ \bU^{(b-1)}({\bbeta}^0) $ could be written as the sum of $ i.i.d. $ random variables, that is, $\bU^{(b-1)}({\bbeta}^0) = \sum_{i \in \mathcal{D}_b^\star} \bu(y_i; \bx_{i}, \bbeta^0).$ Since $ \mathbb{E}[\bU^{(b-1)}({\bbeta}^0)] = 0$ and $  \lVert \bx_i \lVert_\infty \leq K $, according to Hoeffding's inequality, with probability at least $ 1- p^{-3} $, $ \lVert \bU^{(b-1)}({\bbeta}^0) \lVert_\infty \leq \lambda_{b-1}N_{b-1}/8. $
	
	The result of $ |\Delta_4^{(b)}| $ is straightforward according to Theorem 14.5 in \cite{buhlmann2011statistics}.
	
\end{proof}

\bibliographystyle{authordate1}
\bibliography{paper-ref}

\appendix

\end{document}